\documentclass[11pt,reqno]{amsart}
\usepackage[hidelinks]{hyperref}
\usepackage[a4paper,margin=1in]{geometry}
\usepackage[T1]{fontenc}
\usepackage{lmodern}
\usepackage{amsmath,amssymb,amsthm,mathtools,booktabs}
\usepackage{graphicx,pgfplots}
\usepackage{float}
\usepackage[section]{placeins}
\pgfplotsset{compat=1.17}

\numberwithin{equation}{section}
\newcommand{\R}{\mathbb R}
\newcommand{\ud}{\,\mathrm d}
\newcommand{\norm}[1]{\lVert#1\rVert}
\newcommand{\abs}[1]{\lvert#1\rvert}
\newcommand{\St}{S_T}
\newcommand{\Sd}{\mathcal S_\Delta}
\newcommand{\Ac}{\mathcal A}
\newcommand{\Pc}{\mathcal{P}_\Delta}
\newcommand{\Ic}{I_\Delta}

\newtheorem{theorem}{Theorem}[section]
\newtheorem{lemma}[theorem]{Lemma}
\newtheorem{proposition}[theorem]{Proposition}
\newtheorem{corollary}[theorem]{Corollary}
\newtheorem{assumption}[theorem]{Assumption}
\newtheorem{example}[theorem]{Example}
\theoremstyle{remark}
\newtheorem{remark}[theorem]{Remark}

\title[First-order convergence of the minimum action]{First-order convergence of finite-difference approximations
to the minimum action for diffusions with multiplicative noise}\thanks{This work was supported by the National Natural Science Foundation of China (Nos.\ 12671472 and 12201228) and the Fundamental Research Funds for the Central Universities (No.\ 2026BRSXB006).}
\author{Diancong Jin}
\address{School of Mathematics and Statistics, Huazhong University of Science and Technology, Wuhan 430074, China; Hubei Key Laboratory of Engineering Modeling and Scientific Computing, Huazhong University of Science and Technology, Wuhan 430074, China.}
\email{jindc@hust.edu.cn}
\hypersetup{pdfauthor={Diancong Jin}}
\date{}

\begin{document}
\begin{abstract}
We prove first-order convergence of finite-difference approximations to the Freidlin--Wentzell minimum action for diffusions with multiplicative noise, which answers the question of the sharp rate left open in [Hong, Jin, Sheng, \emph{IMA J. Numer. Anal.} 45 (2025), 1501--1536]. As shown in that work, a bound on the action controls only the $L^2$ norm of the velocity, which leads to a half-order estimate for the difference between the discrete and continuous actions. To improve it, we use the discrete Euler--Lagrange equation to derive mesh-independent $L^\infty$ bounds on the velocity of piecewise affine interpolants of stationary discrete paths with bounded action.
For the case of a unique continuous minimizer,
we identify the leading coefficient of the error between the
discrete and continuous minimum actions on uniform grids, and a scalar example shows that the first-order
error bound is sharp. We also prove strong $H^1(0,T;\R^d)$ subsequential convergence of
interpolated discrete minimizers.
Numerical experiments for a two-protein genetic toggle switch confirm first-order convergence of the discrete action for transitions in both directions.
\end{abstract}


\maketitle

% Manuscript sections, tables, and bibliography are inlined below.
% Edit this main file directly for subsequent manuscript revisions.
\section{Introduction}
\label{sec:intro}

Small random perturbations can influence the long-term dynamics of systems. For example, they can drive a dynamical system from one stable state to another. The theory of Freidlin--Wentzell large deviation 
	\cite{FreidlinWentzell,FengKurtz} describes the probabilities and
	pathways of such rare transitions through a variational problem. More precisely, 
consider the small-noise induced diffusion
\begin{equation}
	\ud X_t^\varepsilon
	=b(X_t^\varepsilon)\ud t
	+\sqrt\varepsilon\,\sigma(X_t^\varepsilon)\ud W_t,
	\qquad X_0^\varepsilon=x_0\in\R^d,
	\label{eq:sde}
\end{equation}
where $W$ is a standard $\R^d$-valued Brownian motion,
$\varepsilon>0$ is sufficiently small.
Fix $T>0$ and $x_T\in\R^d$.

For a pointwise invertible diffusion
matrix $\sigma$, the Freidlin--Wentzell action is
\begin{align}
	\St(\phi)&=\frac12\int_0^T
	\abs{\sigma(\phi(t))^{-1}(\dot\phi(t)-b(\phi(t)))}^2\ud t.
	\label{eq:continuous-action}
\end{align}
The large deviation principle for $X^\varepsilon$ in the path space $C([0,T];\R^d)$ indicates that (see \cite{FreidlinWentzell}) for any $\phi\in H^1(0,T;\R^d)$ and $\rho\ll 1$,
\begin{align*}
	\mathbb P\left(\|X^\varepsilon-\phi\|_{L^\infty}<\rho\right)\asymp \exp\left(-\frac{1}{\varepsilon}S_T(\phi)\right),~\text{as}~\varepsilon\to 0.
\end{align*}
Thus, the smaller the action $S_T(\phi)$ is, the more likely $X^\varepsilon$ lies in a tiny tube of $\phi$. 
Particularly, minimizers of $\St$ over $$\Ac:=\{\phi\in H^1(0,T;\R^d):\phi(0)=x_0,\ \phi(T)=x_T\}$$ describe the most likely transition
paths of $X^\varepsilon$ transferring from $x_0$ to $X_T$ in the small-noise limit, and the minimum action
\begin{align}\label{problem1}
m:=\inf_{\phi\in\Ac}\St(\phi)
\end{align}
	determines the exponential rate of the transition probability.
	Minimum actions and minimizing paths have been used, for instance,
	to analyze noise-induced switching between stable gene-expression
	states \cite{RomaEtAl2005}, to choose the modified drift in
	importance sampling so that a rare transition occurs more often in
	simulation \cite{VandenEijndenWeare2012}, and to trace how weak
	random forcing drives two-dimensional Poiseuille flow from its
	laminar state to a coexisting stable travelling wave \cite{WanYuE2015}.

One of the most popular methods, numerically solving the minimum and minimizers of $S_T$, is  the so-called minimum action method.
The minimum action method, first proposed by \cite{E2004}, discretizes the admissible
	paths and minimizes the resulting discrete action. The geometric and
	adaptive variants were later proposed in
	\cite{HeymannVandenEijnden2008,ZhouRenE2008}, a version with optimal
	linear time scaling in \cite{Wan2015}, and an $hp$-adaptive version
	based on a posteriori error estimates in \cite{WanZhengLin2018}.
	Action minimization has also been used for rare-event sampling
	\cite{VandenEijndenWeare2012,GrafkeVandenEijnden2019} and in the
	sharp-interface analysis of the stochastic Allen--Cahn equation
	\cite{KohnOttoReznikoffVandenEijnden2007}.
Instantons can also be computed by solving the associated Euler--Lagrange
equations~\cite{GrafkeGrauerSchafer2015}. Fluctuations around the
	minimizing path then give the prefactors needed for sharp probability
	estimates \cite{SchorleppTongGrafkeStadler2023,SchorleppGrafke2026}
	and for nonequilibrium transition rates \cite{HellerLimmer2024}.
	More recently, transition paths have been parameterized by neural
	networks \cite{Simonnet2023,HanWuGuZhou2026}, an approach that has
	also been used to compute reaction paths and phase diagrams in active
	matter \cite{ZakineSimonnetVandenEijnden2024}. In addition, deep reinforcement
	learning has been used to search path space \cite{LinZhongRen2025},
	and the recent preprint \cite{HuangYueYu2026} proposes a Laguerre
	spectral method on a semi-infinite time interval for computing
	quasi-potentials.

Despite fruitful results on the algorithm design,
rigorous convergence results for these methods are much scarcer. For the fixed-time problem \eqref{problem1}, \cite[Theorem~4.1]{WanYuZhai2018} proved convergence of a conforming finite element approximation with additive noise via $\Gamma$-convergence, and \cite{HongJinSheng2025} obtained error bounds for the minimum action and subsequential convergence of the minimizing paths for a finite-difference method.

This paper continues the analysis of the finite-difference method studied in \cite{HongJinSheng2025}. Let $0=t_0<t_1<\cdots<t_N=T$ be a partition, with $h_n=t_{n+1}-t_n$ and $\Delta=\max_{0\le n<N}h_n$.
The admissible set of discrete paths is
\[
\Ac_\Delta
:=\bigl\{\psi=(\psi_n)_{n=0}^N\in(\R^d)^{N+1}:
\psi_0=x_0,\ \psi_N=x_T\bigr\}.
\]
For $\theta\in[0,1]$ and  \(\psi\in(\mathbb R^d)^{N+1}\), define
\begin{equation}
	\Sd(\psi)=\frac12\sum_{n=0}^{N-1}h_n
	\left|\sigma(\psi_n)^{-1}
	\left(\frac{\psi_{n+1}-\psi_n}{h_n}
	-b((1-\theta)\psi_n+\theta\psi_{n+1})\right)\right|^2.
	\label{eq:intro-discrete}
\end{equation}
The corresponding discrete minimum is
\[
m_\Delta:=\inf_{\psi\in\Ac_\Delta}\Sd(\psi).
\]
For this discretization on nonuniform meshes, \cite[Theorem~3.2]{HongJinSheng2025} obtained an $O(\Delta^{1/2})$ bound for $|m_\Delta-m|$ with multiplicative noise and an $O(\Delta)$ bound with additive noise.
	The numerical experiment for a scalar SDE with multiplicative noise
	in \cite[Section~5.2]{HongJinSheng2025} showed a rate close to one
	for the minimum action, and the sharp convergence order was left open. The main purpose of this paper is to answer the sharp convergence order for $m_\Delta$ through a rigorous analysis.  

The sharp order is of practical interest because errors in the
	action enter small-noise transition probabilities exponentially.
The large deviation principle
\cite[Theorem 2 and Example 1]{ChiariniFischer2014} and the contraction principle give
the logarithmic approximation
\begin{equation}\label{eq:introP}
	\mathbb P\bigl(\abs{X_T^\varepsilon-x_T}<\rho\bigr)
	\asymp \exp(-m/\varepsilon),
\end{equation}
where $\asymp$ denotes logarithmic equivalence, with
$\varepsilon\downarrow0$ followed by $\rho\downarrow0$.
Replacing $m$ by $m_\Delta$ multiplies the exponential factor by
	$\exp(-(m_\Delta-m)/\varepsilon)$, which tends to one as
	$\varepsilon\downarrow0$ provided $|m_\Delta-m|=o(\varepsilon)$.
	With a half-order error bound, this is guaranteed when
	$\Delta=o(\varepsilon^2)$; with a first-order bound, $\Delta=o(\varepsilon)$
	suffices. The sharp order therefore indicates how small the time step
	needs to be relative to the noise level.

The main difficulty in proving first-order convergence of the
minimum action is that, for a general admissible discrete path
$\psi$ with $\Sd(\psi)\le M$, the action bound provides only
the displacement estimate
\[
|\psi_{n+1}-\psi_n|\le C_Mh_n^{1/2}
\]
(see Lemma~\ref{lem:discrete-sublevel}).
For state-dependent $\sigma$, estimating the variation of
$\sigma^{-1}$ along each interpolated segment from this bound
gives only an estimate of order $h_n^{1/2}$, as shown in \cite{HongJinSheng2025}.
To improve this estimate, we use the discrete Euler--Lagrange
	equation satisfied by stationary paths.
For a path $\Psi$ that is stationary for $\Sd$ on $\Ac_\Delta$ and satisfies $\Sd(\Psi)\le M$, we derive the discrete Euler--Lagrange equation and apply it to obtaining
\[
|\Psi_{n+1}-\Psi_n|\le C_M h_n
\]
(see Proposition~\ref{prop:discrete-reg}).
With this uniform bound, we further prove that $\Sd(\Psi)$ and $\St(\Ic\Psi)$ differ
	by $O(\Delta)$ in Lemma \ref{lem:affine-consistency}, where $\Ic\Psi$ is the continuous piecewise affine
	interpolant of $\Psi$. Applying Lemma \ref{lem:affine-consistency} to discrete
	minimizers, together with a sampling estimate for continuous
	minimizers in Lemma \ref{lem:sampling}, gives the main result (Theorem~\ref{thm:main})
\[
|m_\Delta-m|\le C\Delta.
\]

Building on Theorem~\ref{thm:main}, we also prove that interpolated
discrete minimizers have subsequences converging strongly in
$H^1(0,T;\R^d)$ to minimizers of $\St$ over $\Ac$
(Proposition~\ref{prop:minimizer-convergence}).
This strengthens the weak $H^1$ subsequential convergence
established in \cite[Theorem~3.3]{HongJinSheng2025}.
If, in addition, $b$ and $\sigma$ are of class $C^2$ and
$\St$ has a unique minimizer $\phi^*$ over $\Ac$, then
\[
m_\Delta-m
=\Delta\mathcal B_\theta(\phi^*)+o(\Delta)
\]
on uniform grids with fixed $\theta$
(Theorem~\ref{thm:bias}).
The coefficient $\mathcal B_\theta$, given in
	\eqref{eq:bias-functional}, has two parts: one comes from
	evaluating the drift at the $\theta$-weighted nodes; the other, from
	evaluating the diffusion matrix at the left endpoints.
	Example~\ref{ex:sharp} gives a scalar diffusion with
	$\mathcal B_\theta(\phi^*)\ne0$, so the first-order estimate for
	the minimum action is sharp in general. 
	On uniform grids, if \(A=(\sigma\sigma^\top)^{-1}\) is constant and the drift is evaluated at the midpoint, then \(\mathcal B_{1/2}\equiv0\). Lemma 4.6 and (4.25) then give \(|m_\Delta-m|\le C\Delta^2\), as stated in Corollary \ref{cor:midpoint-minimum-action}.   
If $\St$ also satisfies the local quadratic growth condition
\eqref{eq:path-quadratic-growth} at $\phi^*$, we prove
\[
\norm{\Ic\psi^\Delta-\phi^*}_{H^1}\le C\Delta
\]
for every global minimizer $\psi^\Delta$ of $\Sd$ over
$\Ac_\Delta$ on sufficiently fine uniform grids
(Theorem~\ref{thm:path-rate}).

The paper is organized as follows.
Section~\ref{sec:framework} states the assumptions on $b$ and
$\sigma$, records bounds for continuous and discrete paths with
bounded action, and recalls the existence of minimizers.
In Section~\ref{sec:regularity},  we derive the continuous and discrete
Euler--Lagrange equations and prove uniform regularity estimates
for stationary paths with bounded action.
Section~\ref{sec:convergence} proves first-order convergence of
the minimum action and strong $H^1$ subsequential convergence of
interpolated discrete minimizers. It also studies the leading bias
on uniform grids, the sharpness of the first-order minimum-action
estimate, and conditions for first-order convergence of the
minimizing paths. Finally, in
Section~\ref{sec:numerics}, we present a two-protein genetic toggle
model and report the computed switching paths and the action
and path errors for transitions in both directions.

\section{Preliminaries}
\label{sec:framework}
We begin with some notations. 
We write $|\cdot|$ for the absolute value of a scalar and the
Euclidean norm of a vector, and $\norm{\cdot}$ for the induced operator
norm on matrices and linear maps.
For a Banach space $X$, we write $u_j\rightharpoonup u$ in $X$
for weak convergence and $u_j\to u$ in $X$ for convergence in norm.
The generic finite constant $C$ may
vary between occurrences.
Dependence on fixed data such as $T,x_0,x_T,b,\sigma$ is suppressed.
We use subscripts to indicate additional dependence, as in $C_M$ or $C_{R,K}$.
Unless stated otherwise, all constants are independent of the partition,
$N$, and $\theta\in[0,1]$.
For $R>0$, we write
$\overline{B_R(0)}:=\{x\in\R^d:\abs{x}\le R\}$ for the closed Euclidean
ball centered at the origin. Balls in other spaces will be specified
explicitly. In addition, 
$H_0^1(0,T;\R^d):=\{\eta\in H^1(0,T;\R^d):\eta(0)=\eta(T)=0\}$ denotes the space of endpoint-preserving variations.

Recall the partition $0=t_0<t_1<\cdots<t_N=T$, with
$h_n=t_{n+1}-t_n$ and $\Delta=\max_{0\le n<N}h_n$.
For a discrete path $\psi=(\psi_n)_{n=0}^N\in(\R^d)^{N+1}$ and
$\theta\in[0,1]$, define the difference quotient and the weighted node by
\begin{equation}
(\delta\psi)_n=\frac{\psi_{n+1}-\psi_n}{h_n},
\qquad
\psi_n^\theta=(1-\theta)\psi_n+\theta\psi_{n+1},
\qquad 0\le n<N.
\label{eq:nodal-notation}
\end{equation}
The continuous piecewise affine interpolation operator $\Ic$ is defined by
\[
(\Ic\psi)(t):=\psi_n+(t-t_n)(\delta\psi)_n,
\qquad t\in[t_n,t_{n+1}],\quad 0\le n<N.
\]
For $\phi\in C([0,T];\R^d)$, we define the nodal sampling operator by
\begin{equation}\label{eq:P}
\Pc\phi:=(\phi(t_n))_{n=0}^N.
\end{equation}
The interpolant $\Ic\psi$ belongs to $H^1(0,T;\R^d)$.
It is easy to see that both operators preserve the prescribed endpoint values, so
$\Ic(\Ac_\Delta)\subset\Ac$ and $\Pc(\Ac)\subset\Ac_\Delta$.

\begin{assumption}
\label{ass:main}
The coefficients $b\in C^1(\R^d;\R^d)$ and
$\sigma\in C^1(\R^d;\R^{d\times d})$ 
have at most linear growth: there exists $L_0\ge1$ such that
\begin{equation}
\abs{b(x)}+\norm{\sigma(x)}\le L_0(1+\abs{x}),\qquad x\in\R^d.
\label{eq:growth}
\end{equation}
Moreover, $\det\sigma(x)\ne0$ for every $x\in\R^d$.
\end{assumption}

For example, the scalar
diffusion $\sigma(x)=(1+x^2)^{-1/2}$ is admissible, although
$\inf_{x\in\R}\sigma(x)=0$.
Set $\Delta_0:=(2L_0)^{-1}$. Under Assumption \ref{ass:main}, $b$  and $\sigma$ are locally Lipschitz, which means that for every $R>0$, there exists a constant $C_R>0$ such that
\[
\abs{b(x)-b(y)}+\norm{\sigma(x)-\sigma(y)}
\le C_R\abs{x-y},\qquad \abs{x},\abs{y}\le R.
\]
Pointwise invertibility and continuity of $\sigma$ imply
\[
\sup_{\abs{x}\le R}\norm{\sigma^{-1}(x)}<\infty,
\qquad R>0.
\]
The inverse $\sigma^{-1}$ is locally Lipschitz as well.
%\begin{align*}
%\norm{\sigma^{-1}(x)-\sigma^{-1}(y)}
%&\le C_R\bigg(\sup_{\abs{z}\le R}\norm{\sigma^{-1}(z)}\bigg)^2\abs{x-y}.
%\end{align*}
%Thus $\sigma^{-1}$ is locally Lipschitz as well.

With this notation, the finite-difference action takes the form
\begin{equation*}
\Sd(\psi)=\frac12\sum_{n=0}^{N-1}h_n
\left|\sigma(\psi_n)^{-1}
\bigl((\delta\psi)_n-b(\psi_n^\theta)\bigr)\right|^2.
\end{equation*}
The following two lemmas provide continuous and discrete sublevel bounds.

\begin{lemma}
\label{lem:continuous-sublevel}
Under Assumption~\ref{ass:main}, for every $M\ge0$ there exists $C_M$
such that every $\phi\in H^1(0,T;\R^d)$ with $\phi(0)=x_0$ and
$\St(\phi)\le M$ satisfies
\begin{equation}
\norm{\phi}_{L^\infty}+\norm{\dot\phi}_{L^2}\le C_M.
\label{eq:continuous-sublevel}
\end{equation}
\end{lemma}

\begin{proof}
The proof of \cite[Proposition~2.3]{HongJinSheng2025} uses global
Lipschitz continuity only through the linear growth bound
\cite[(2.2)]{HongJinSheng2025}, which is \eqref{eq:growth} here.
Its $H^1(0,T;\R^d)$ estimate and
\cite[Proposition~2.2(2)]{HongJinSheng2025} imply
\eqref{eq:continuous-sublevel}.
\end{proof}

The discrete estimate retain the original cell lengths and is uniform
in $\theta\in[0,1]$.

\begin{lemma}
\label{lem:discrete-sublevel}
Under Assumption~\ref{ass:main}, let $\Delta\le\Delta_0$.
For every $M\ge0$ there exists $C_M$ such that every discrete path
$\psi$ with $\psi_0=x_0$ and $\Sd(\psi)\le M$ satisfies
\begin{equation}
\max_{0\le n\le N}\abs{\psi_n}
+\left(\sum_{n=0}^{N-1}h_n\abs{(\delta\psi)_n}^2\right)^{1/2}
\le C_M.
\label{eq:discrete-sublevel}
\end{equation}
\end{lemma}

\begin{proof}
On the interpolant $\Ic\psi$, the action
in \cite[Lemma~2.8]{HongJinSheng2025} equals $\Sd(\psi)$ by
\cite[(2.9)]{HongJinSheng2025}. The proof of that lemma also uses only
the linear growth bound \eqref{eq:growth}. With $L=L_0$, it bounds the
nodes and the derivative for $\Delta\le(2L_0)^{-1}$ and yields
\eqref{eq:discrete-sublevel}, uniformly in the partition and
$\theta\in[0,1]$.
\end{proof}

We next establish the existence of minimizers and a common action bound.
\begin{proposition}
\label{prop:existence}
Under Assumption~\ref{ass:main}, $\St$ attains its minimum over $\Ac$.
For every partition with $\Delta\le\Delta_0$, $\Sd$ attains its minimum
over $\Ac_\Delta$, and
\begin{equation}
\max\{m,m_\Delta\}\le M_0<\infty.
\label{eq:comparison-level}
\end{equation}
\end{proposition}

\begin{proof}
Set $\ell(t):=x_0+tv$ for $t\in[0,T]$, where
\(
v:=(x_T-x_0)/T.
\)
For $0\le n<N$, we have $
\dot\ell=(\delta\Pc\ell)_n=v$ and
$
(\Pc\ell)_n^\theta=\ell(t_n+\theta h_n)\in \ell([0,T])$.
Since $b$ and $\sigma^{-1}$ are bounded on $\ell([0,T])$,
\[
\begin{aligned}
\max\{m,m_\Delta\}
\le\max\{\St(\ell),\Sd(\Pc\ell)\}\le CT(1+\abs{v}^2)=:M_0.
\end{aligned}
\]

The existence follows by the direct method, as in
\cite[Lemmas~2.5 and~2.11]{HongJinSheng2025}. We
choose a minimizing sequence $(\phi_j)\subset\Ac$ with
$\St(\phi_j)\le M_0+1$. Then Lemma~\ref{lem:continuous-sublevel} and the
compact embedding of $H^1(0,T;\R^d)$ into $C([0,T];\R^d)$ yield,
after passing to a subsequence,
\[
\phi_j\rightharpoonup\phi\quad\text{in }H^1(0,T;\R^d),
\qquad
\phi_j\to\phi\quad\text{in }C([0,T];\R^d),
\qquad
\phi\in\Ac.
\]
The uniform bound gives
$\phi_j([0,T])\cup\phi([0,T])\subset\mathcal K:=\overline{B_R(0)}\subset\R^d$
for some $R>0$ independent of $j$. It follows from the uniform continuity of $b$ and
$\sigma^{-1}$ on this compact set and the preceding
convergences that
\[
\begin{aligned}
\sigma(\phi_j)^{-1}\bigl(\dot\phi_j-b(\phi_j)\bigr)\rightharpoonup
\sigma(\phi)^{-1}\bigl(\dot\phi-b(\phi)\bigr)
\quad\text{in }L^2(0,T;\R^d).
\end{aligned}
\]
Hence, the weak lower semicontinuity of the $L^2$ norm gives
\[
m\le\St(\phi)\le\liminf_{j\to\infty}\St(\phi_j)=m.
\]

For the discrete problem, fix a partition with $\Delta\le\Delta_0$
and choose a minimizing sequence $(\psi^{(j)})\subset\Ac_\Delta$
such that
\[
\Sd(\psi^{(j)})\rightarrow m_\Delta,
\qquad
\Sd(\psi^{(j)})\le M_0+1.
\]
For the fixed partition, Lemma~\ref{lem:discrete-sublevel} bounds
$(\psi^{(j)})$ in $(\R^d)^{N+1}$. Accordingly, there exists a subsequence, still denoted by
$(\psi^{(j)})$, such that
\[
\psi_n^{(j)}\rightarrow\psi_n,
\qquad 0\le n\le N.
\]
The endpoint conditions are preserved when passing to the limit, and thus 
$\psi:=(\psi_0,\ldots,\psi_N)\in\Ac_\Delta$.
Moreover, since each $h_n>0$ is fixed,
\[
(\delta\psi^{(j)})_n\rightarrow(\delta\psi)_n,
\qquad
(\psi^{(j)})_n^\theta\rightarrow\psi_n^\theta,
\qquad 0\le n<N.
\]
Continuity of $b$ and $\sigma^{-1}$ yields
\(
\Sd(\psi)
=\lim_{j\to\infty}\Sd(\psi^{(j)})
=m_\Delta.
\)
Thus $\psi$ is a discrete minimizer.
\end{proof}

For comparison, we recall \cite[Theorem~3.2]{HongJinSheng2025} under
its global Lipschitz assumptions. The continuous-path discrete action
used there has the same minimum $m_\Delta$ as $\Sd$ by
\cite[Lemma~2.6]{HongJinSheng2025}.

\begin{proposition}
\label{prop:previous-rates}
Suppose, in addition to Assumption~\ref{ass:main}, that $b$ and
$\sigma$ are globally Lipschitz. For every sufficiently fine partition
and every $\theta\in[0,1]$,
\begin{equation}
|m_\Delta-m|\le C\Delta^{1/2}.
\label{eq:previous-half-order}
\end{equation}
If $\sigma$ is an invertible constant matrix, then
\begin{equation}
|m_\Delta-m|\le C\Delta.
\label{eq:previous-additive-order}
\end{equation}
\end{proposition}

In Theorem~\ref{thm:main} we will use the Euler--Lagrange equations of minimizing
paths to prove a first-order estimate for state-dependent diffusion under \(C^1\) regularity and linear growth, without requiring global Lipschitz continuity.

% ===== Section 3: Regularity of continuous and discrete paths =====
\section{Regularity of continuous and discrete paths}
\label{sec:regularity}

Define the matrix $A(x)=(\sigma(x)\sigma(x)^{\top})^{-1}$ and set
\begin{equation}
F(x,v)=\tfrac12(v-b(x))^{\top}A(x)(v-b(x)),\quad x,v\in\R^d.
\label{eq:AF}
\end{equation}
Then we rewrite $S_T$ as
\begin{equation}
\St(\phi)=\int_0^T F(\phi,\dot\phi)\ud t.
\label{eq:action-integral}
\end{equation}
Under Assumption~\ref{ass:main}, for every $R>0$ there exist
$0<c_R\le C_R<\infty$ such that, for $\abs{x}\le R$ and $\xi\in\R^d$,
\begin{gather}
c_R\abs{\xi}^2
\le \xi^{\top}A(x)\xi\le C_R\abs{\xi}^2,
\label{eq:local-elliptic}\\
\norm{DA(x)}+\norm{D(A^{-1})(x)}+\abs{b(x)}+\norm{Db(x)}
\le C_R.
\label{eq:local-coefficient-bounds}
\end{gather}
By Lemmas~\ref{lem:continuous-sublevel}
and~\ref{lem:discrete-sublevel}, continuous paths with $\phi(0)=x_0$
and $\St(\phi)\le M$, and discrete paths with $\psi_0=x_0$,
$\Sd(\psi)\le M$, and $\Delta\le\Delta_0$, take values in a common
closed ball $\overline{B_{R_M}(0)}\subset\R^d$. The radius is independent
of the path, the partition, and $\theta$. By convexity, the piecewise affine
interpolants of these discrete paths also take values in this closed ball.

The next proposition states the regularity and uniform bounds for
continuous stationary paths with bounded action.

\begin{proposition}
\label{prop:continuous-reg}
Under Assumption~\ref{ass:main}, let $\Phi\in\Ac$ satisfy
\[
\St(\Phi)\le M,
\qquad
D\St(\Phi)[\eta]=0
\quad\text{for all }\eta\in H_0^1(0,T;\R^d).
\]
Then $\Phi\in C^2([0,T];\R^d)$, and the momentum
\(
P:=\nabla_vF(\Phi,\dot\Phi)
\)
belongs to $C^1([0,T];\R^d)$. Moreover,
\begin{equation}
\norm{P}_{L^\infty}+\norm{\dot P}_{L^\infty}
+\norm{\dot\Phi}_{L^\infty}+\norm{\ddot\Phi}_{L^\infty}
\le C_M.
\label{eq:continuous-reg}
\end{equation}
\end{proposition}

\begin{proof}
For $\eta\in C_c^\infty((0,T);\R^d)$ and almost every $t\in(0,T)$,
the chain rule reads
\begin{align*}
\frac{\mathrm d}{\mathrm d\varepsilon}
F(\Phi+\varepsilon\eta,\dot\Phi+\varepsilon\dot\eta)=
\nabla_xF(\Phi+\varepsilon\eta,\dot\Phi+\varepsilon\dot\eta)\cdot\eta+
\nabla_vF(\Phi+\varepsilon\eta,\dot\Phi+\varepsilon\dot\eta)\cdot\dot\eta.
\end{align*}
By the embedding
$H^1(0,T;\R^d)\hookrightarrow C([0,T];\R^d)$, all values
$\Phi(t)+\varepsilon\eta(t)$, with $t\in[0,T]$ and $|\varepsilon|\le1$,
lie in the closed ball $\overline{B_R(0)}\subset\R^d$, where
$R:=1+\norm{\Phi}_{L^\infty}+\norm{\eta}_{L^\infty}$.
This ball is compact and independent of $t$ and $\varepsilon$. Then on this ball, $A,DA,b,Db$ are bounded by
\eqref{eq:local-elliptic} and \eqref{eq:local-coefficient-bounds}.
Since $F(x,v)$ is quadratic in $v$ and both $\eta$ and $\dot\eta$ are
bounded, one has
\[
\sup_{|\varepsilon|\le1}
\left|
\frac{\mathrm d}{\mathrm d\varepsilon}
F(\Phi+\varepsilon\eta,\dot\Phi+\varepsilon\dot\eta)
\right|
\le C_{\Phi,\eta}(1+\abs{\dot\Phi}^2)\in L^1(0,T;\R).
\]
It follows from the stationarity of $\Phi$
that for every $\eta\in C_c^\infty((0,T);\R^d)$,
\begin{align}
0=D\St(\Phi)[\eta]
=\int_0^T
\left.\frac{\mathrm d}{\mathrm d\varepsilon}
F(\Phi+\varepsilon\eta,\dot\Phi+\varepsilon\dot\eta)
\right|_{\varepsilon=0}\ud t=\int_0^T\bigl(P\cdot\dot\eta
+\nabla_xF(\Phi,\dot\Phi)\cdot\eta\bigr)\ud t.
\label{eq:first-variation}
\end{align}
Thus $\dot P=\nabla_xF(\Phi,\dot\Phi)$ in distributions.

By Lemma~\ref{lem:continuous-sublevel}, 
$\norm{\Phi}_{L^\infty}+\norm{\dot\Phi}_{L^2}\le C_M$.
For $r=v-b(x)$ and $\xi\in\R^d$, it holds that
\begin{align}
\nabla_vF(x,v)&=A(x)r, \notag\\
\nabla_xF(x,v)\cdot\xi
&=\tfrac12r^{\top}DA(x)[\xi]r
-(A(x)r)\cdot Db(x)\xi.
\label{eq:metric-derivative}
\end{align}
In particular, $P=A(\Phi)\bigl(\dot\Phi-b(\Phi)\bigr)$.
Applying \eqref{eq:local-elliptic} and
\eqref{eq:local-coefficient-bounds} to these expressions, we obtain
\begin{equation}
\norm{P}_{L^2}
+\norm{\nabla_xF(\Phi,\dot\Phi)}_{L^1}\le C_M.
\label{eq:continuous-momentum-integrability}
\end{equation}
Since $P\in L^2(0,T;\R^d)$ and its weak derivative
$\nabla_xF(\Phi,\dot\Phi)\in L^1(0,T;\R^d)$,
we have $P\in W^{1,1}(0,T;\R^d)$. Accordingly, the Sobolev embedding
\(
W^{1,1}(0,T;\R^d)\hookrightarrow C([0,T];\R^d)
\)
and \eqref{eq:continuous-momentum-integrability} yield
\begin{align}
\norm{P}_{L^\infty}
\le C_T\bigl(\norm{P}_{L^1}+\norm{\dot P}_{L^1}\bigr)\le C_T\bigl(\norm{P}_{L^2}
+\norm{\nabla_xF(\Phi,\dot\Phi)}_{L^1}\bigr)
\le C_M.
\label{eq:continuous-TV}
\end{align}

Note that $\dot\Phi=b(\Phi)+A(\Phi)^{-1}P$ almost everywhere.
Taking the essential supremum and using \eqref{eq:local-elliptic},
\eqref{eq:local-coefficient-bounds}, and \eqref{eq:continuous-TV}, we obtain
\begin{equation}
\norm{\dot\Phi}_{L^\infty}
\le \norm{b(\Phi)}_{L^\infty}
+\norm{A(\Phi)^{-1}}_{L^\infty}\norm{P}_{L^\infty}
\le C_M.
\label{eq:continuous-velocity-bound}
\end{equation}
Since $\Phi\in H^1(0,T;\R^d)$ also has a continuous representative,
the function $b(\Phi)+A(\Phi)^{-1}P$ is continuous.
In integral form, the identity for $\dot\Phi$ reads
$$
\Phi(t)=x_0+\int_0^t
\bigl(b(\Phi(\tau))+A(\Phi(\tau))^{-1}P(\tau)\bigr)\ud\tau,
\qquad t\in[0,T].
$$
It follows that $\Phi\in C^1([0,T];\R^d)$.
As $F\in C^1(\R^d\times\R^d;\R)$, the function $\nabla_xF(\Phi,\dot\Phi)$ is now
continuous. The identity $\dot P=\nabla_xF(\Phi,\dot\Phi)$ almost
everywhere then implies
$P\in C^1([0,T];\R^d)$.
Thus $P$ satisfies the continuous Euler--Lagrange equation for $\St$, namely,
\begin{equation}
\dot P=\frac{\mathrm d}{\mathrm dt}\nabla_vF(\Phi,\dot\Phi)
=\nabla_xF(\Phi,\dot\Phi),\qquad 0<t<T.
\label{eq:continuous-EL}
\end{equation}
Next, \eqref{eq:continuous-EL} and \eqref{eq:metric-derivative}, together
with \eqref{eq:continuous-velocity-bound}, imply
\begin{equation}
\norm{\dot P}_{L^\infty}
=\norm{\nabla_xF(\Phi,\dot\Phi)}_{L^\infty}
\le C_M\bigl(1+\norm{\dot\Phi}_{L^\infty}^2\bigr)
\le C_M.
\label{eq:continuous-momentum-derivative-bound}
\end{equation}
Since $b\in C^1(\R^d;\R^d)$ and
$A^{-1}\in C^1(\R^d;\R^{d\times d})$, it follows that
$\Phi\in C^2([0,T];\R^d)$, with
\[
\ddot\Phi
=Db(\Phi)\dot\Phi
+D(A^{-1})(\Phi)[\dot\Phi]P
+A(\Phi)^{-1}\dot P.
\]
Estimating these terms by \eqref{eq:continuous-TV},
\eqref{eq:continuous-velocity-bound}, \eqref{eq:continuous-momentum-derivative-bound},
and \eqref{eq:local-coefficient-bounds} proves \eqref{eq:continuous-reg}.
\end{proof}
Note taht every minimizer $\phi^*$ of $\St$ over $\Ac$ is stationary under
endpoint-preserving variations and satisfies $\St(\phi^*)\le M_0$
by Proposition~\ref{prop:existence}.
Applying Proposition~\ref{prop:continuous-reg} with $\Phi=\phi^*$
and $M=M_0$ improves the $H^1(0,T;\R^d)$ bound in
Lemma~\ref{lem:continuous-sublevel} to $C^2([0,T];\R^d)$ regularity
with uniform bounds on the first two derivatives.

We next derive analogous bounds on the difference quotients and discrete momenta
of paths that are stationary for $\Sd$ on $\Ac_\Delta$.
For $\psi=(\psi_n)_{n=0}^N\in\Ac_\Delta$, set
\begin{equation}
r_n:=(\delta\psi)_n-b(\psi_n^\theta),
\qquad
p_n:=A(\psi_n)r_n,
\qquad 0\le n<N.
\label{eq:discrete-momentum}
\end{equation}
With this notation, the discrete action reads
\begin{equation}
\Sd(\psi)=\frac12\sum_{n=0}^{N-1}h_n
r_n^{\top}A(\psi_n)r_n,
\qquad \psi\in\Ac_\Delta.
\label{eq:discrete-action-quadratic}
\end{equation}
The next lemma records how $\Sd$ changes under a variation of a single
interior node. 
\begin{lemma}
\label{lem:discrete-EL}
Let $\psi\in\Ac_\Delta$. Then, for every $1\le n<N$ and $\xi\in\R^d$,
\begin{align}
(p_n-p_{n-1})\cdot\xi
={}&\tfrac12h_nr_n^{\top}DA(\psi_n)[\xi]r_n \notag\\
&-h_n(1-\theta)p_n\cdot Db(\psi_n^\theta)\xi \notag\\
&-h_{n-1}\theta p_{n-1}\cdot
Db(\psi_{n-1}^\theta)\xi-D_{\psi_n}\Sd(\psi)[\xi].
\label{eq:discrete-EL}
\end{align}
\end{lemma}

\begin{proof}
Let $1\le n<N$ and $\xi\in\R^d$. The derivative of
\eqref{eq:discrete-action-quadratic} with respect to $\psi_n$ in
the direction $\xi$ is
\begin{align}
D_{\psi_n}\Sd(\psi)[\xi]
%&=\sum_{i=0}^{N-1}h_i
%\bigl(D_{\psi_n}r_i[\xi]\bigr)^{\top}A(\psi_i)r_i\\
%&\quad+\frac12\sum_{i=0}^{N-1}h_i
%r_i^{\top}D_{\psi_n}\bigl(A(\psi_i)\bigr)[\xi]r_i\\
&=h_np_n\cdot D_{\psi_n}r_n[\xi]
+h_{n-1}p_{n-1}\cdot D_{\psi_n}r_{n-1}[\xi]+\frac12h_n r_n^{\top}DA(\psi_n)[\xi]r_n.
\label{eq:discrete-node-variation}
\end{align}
Since $r_n=(\delta\psi)_n-b(\psi_n^\theta)$, we obtain
\begin{align*}
D_{\psi_n}r_n[\xi]
&=-h_n^{-1}\xi-(1-\theta)Db(\psi_n^\theta)\xi,\\
D_{\psi_n}r_{n-1}[\xi]
&=h_{n-1}^{-1}\xi-\theta Db(\psi_{n-1}^\theta)\xi.
\end{align*}
Substituting both derivatives into \eqref{eq:discrete-node-variation}
and collecting the momentum terms proves \eqref{eq:discrete-EL}.
\end{proof}

In particular, $\psi$ is stationary for $\Sd$ on $\Ac_\Delta$ if and only if
$\nabla_{\psi_n}\Sd(\psi)=0$ for all $1\le n<N$. In this case,
\eqref{eq:discrete-EL} is the discrete Euler--Lagrange equation.
We next show that discrete paths with bounded action and controlled interior gradients
satisfy the following uniform bounds.

\begin{lemma}
\label{lem:residual}
Under Assumption~\ref{ass:main}, let
$\Delta\le\Delta_0=(2L_0)^{-1}$ and let
$\psi\in\Ac_\Delta$ satisfy
$\Sd(\psi)\le M$. At each interior node $1\le n<N$,
\begin{equation}
\abs{p_n-p_{n-1}}
\le \abs{\nabla_{\psi_n}\Sd(\psi)}
+C_M\bigl(h_n\abs{p_n}^2+h_n\abs{p_n}+h_{n-1}\abs{p_{n-1}}\bigr).
\label{eq:nodal-residual-increment}
\end{equation}
If $\sum_{n=1}^{N-1}\abs{\nabla_{\psi_n}\Sd(\psi)}\le R$, then
$$
\max_{0\le n<N}\abs{p_n}
+\sum_{n=1}^{N-1}\abs{p_n-p_{n-1}}
+\max_{0\le n<N}\abs{(\delta\psi)_n}
\le C_M(1+R).
$$
\end{lemma}

\begin{proof}
By Lemma~\ref{lem:discrete-sublevel}, the nodes
$\psi_n$ and their convex combinations $\psi_n^\theta$ lie in a common
closed ball $\overline{B_{R_M}(0)}\subset\R^d$. In the
first-variation identity \eqref{eq:discrete-EL}, substitute
$r_n=A(\psi_n)^{-1}p_n$ and use the local bounds on $DA$, $A^{-1}$ and
$Db$ in \eqref{eq:local-elliptic} and
\eqref{eq:local-coefficient-bounds}. Taking the supremum over $\abs{\xi}=1$
then yields \eqref{eq:nodal-residual-increment}.

By \eqref{eq:local-elliptic}, all  eigenvalues of the symmetric positive definite matrix
$A(\psi_n)$  lie in $(0,C_M]$.
Each such eigenvalue satisfies $\lambda^2\le C_M\lambda$, so 
$
A(\psi_n)^2\le C_MA(\psi_n)
$
in the sense of quadratic forms. Recalling $p_n=A(\psi_n)r_n$, we obtain
$$
\abs{p_n}^2
=r_n^{\top}A(\psi_n)^2r_n
\le C_Mr_n^{\top}A(\psi_n)r_n.
$$
Multiplying by $h_n$, summing over the cells, and using
\eqref{eq:discrete-action-quadratic} yields
\begin{align}
E_\Delta
:=\sum_{n=0}^{N-1}h_n\abs{p_n}^2
&\le C_M\sum_{n=0}^{N-1}
h_nr_n^{\top}A(\psi_n)r_n =2C_M\Sd(\psi)\le C_M.
\label{eq:momentum-energy}
\end{align}
Now suppose $\sum_{n=1}^{N-1}\abs{\nabla_{\psi_n}\Sd(\psi)}\le R$.
Summing \eqref{eq:nodal-residual-increment} over the interior nodes, using the Cauchy--Schwarz inequality and \eqref{eq:momentum-energy}, we conclude
\begin{align}
\sum_{n=1}^{N-1}\abs{p_n-p_{n-1}}
&\le R+C_M\sum_{n=0}^{N-1}h_n\bigl(\abs{p_n}^2+\abs{p_n}\bigr) \le R+C_M\bigl(E_\Delta+\sqrt{T E_\Delta}\bigr)
\le R+C_M.
\label{eq:momentum-total-variation}
\end{align}
Since $\sum_{n=0}^{N-1}h_n=T$, there exists
$j\in\{0,\ldots,N-1\}$ such that
\begin{equation}
\abs{p_j}^2
\le\frac1T\sum_{n=0}^{N-1}h_n\abs{p_n}^2
=\frac{E_\Delta}{T}.
\label{eq:momentum-average-node}
\end{equation}
Combining \eqref{eq:momentum-total-variation} and
\eqref{eq:momentum-average-node} yields
\begin{align}
\max_{0\le n<N}\abs{p_n}
\le \abs{p_j}+\sum_{k=1}^{N-1}\abs{p_k-p_{k-1}} \le\left(\frac{E_\Delta}{T}\right)^{1/2}+R+C_M.
\label{eq:energy-TV}
\end{align}
Finally, $(\delta\psi)_n=b(\psi_n^\theta)+A(\psi_n)^{-1}p_n$.
Using the bounds on $b$ and $A^{-1}$ in \eqref{eq:local-coefficient-bounds}
and \eqref{eq:local-elliptic}, respectively, we conclude that
$\max_n\abs{(\delta\psi)_n}\le C_M(1+R)$.
\end{proof}

The following proposition shows that stationary paths for $\Sd$ on $\Ac_\Delta$ satisfy mesh-independent
bounds on their momenta and difference quotients, which is the discrete counterpart of
Proposition~\ref{prop:continuous-reg}.

\begin{proposition}
\label{prop:discrete-reg}
Under Assumption~\ref{ass:main}, let $\Delta\le\Delta_0$.
If $\Psi\in\Ac_\Delta$ is stationary for $\Sd$ and
$\Sd(\Psi)\le M$, then
\begin{equation}
\max_{0\le n<N}\abs{P_n}
+\sum_{n=1}^{N-1}\abs{P_n-P_{n-1}}
+\max_{0\le n<N}\abs{(\delta\Psi)_n}
\le C_M,
\label{eq:discrete-reg}
\end{equation}
where $
P_n:=A(\Psi_n)\bigl((\delta\Psi)_n-b(\Psi_n^\theta)\bigr)$, $0\le n<N$. 
Moreover,
\begin{align}
\abs{\Psi_{n+1}-\Psi_n}
&\le C_Mh_n,
&&0\le n<N, \notag\\
\abs{P_n-P_{n-1}}
&\le C_M(h_n+h_{n-1}),
&&1\le n<N.
\label{eq:increments}
\end{align}
\end{proposition}

\begin{proof}
We take $\psi=\Psi$ in \eqref{eq:discrete-momentum} and Lemma~\ref{lem:residual}. Then we have $p_n=P_n$ and $\nabla_{\psi_n}\mathcal S_\Delta(\Psi)=0,~ 1\le n<N$. With $R=0$, the
estimate in the lemma reduces to \eqref{eq:discrete-reg}.
Multiplying $\abs{(\delta\Psi)_n}\le C_M$ by $h_n$ yields the first bound in
\eqref{eq:increments}. Inserting $\abs{P_k}\le C_M$ in
\eqref{eq:nodal-residual-increment} proves the second bound in
\eqref{eq:increments}.
\end{proof}

% ===== Section 4: Convergence of the minimum action and minimizing paths =====
\section{Convergence of the minimum action and minimizing paths}
\label{sec:convergence}

In this section, we establish first-order convergence of the minimum action
and strong subsequential convergence of interpolated discrete minimizers.
For uniform grids, we also study the leading bias of $m_\Delta-m$ and the convergence
rate of minimizing paths.

\subsection{Convergence and error estimates}
\label{sec:basic-convergence}
\leavevmode\par

Under Assumption~\ref{ass:main}, we first compare the actions on
arbitrary time partitions with $\theta\in[0,1]$, using the interpolation
operator $\Ic$ defined in Section~\ref{sec:framework}.
The slope bounds in
Propositions~\ref{prop:continuous-reg} and~\ref{prop:discrete-reg}
allow us to apply the following estimate to minimizing paths.

\begin{lemma}
\label{lem:affine-consistency}
For every $R,K\ge0$, let $\psi\in\Ac_\Delta$ satisfy
\[
\max_{0\le n\le N}\abs{\psi_n}\le R,
\qquad
\max_{0\le n<N}\abs{(\delta\psi)_n}\le K.
\]
Then
\begin{equation}
|\St(\Ic\psi)-\Sd(\psi)|
\le C_{R,K}T\Delta.
\label{eq:affine-consistency}
\end{equation}
\end{lemma}

\begin{proof}
Since $(\Ic\psi)'(t)=(\delta\psi)_n$ for $t\in(t_n,t_{n+1})$, it follows from
\eqref{eq:action-integral} and \eqref{eq:discrete-action-quadratic} that 
\begin{align}\label{eq:Slcphi}
|\St(\Ic\psi)-\Sd(\psi)|
&=\biggl|\sum_{n=0}^{N-1}\int_{t_n}^{t_{n+1}}
\Bigl(F((\Ic\psi)(t),(\delta\psi)_n)
-\tfrac12 r_n^{\top}A(\psi_n)r_n\Bigr)\ud t\biggr|.
\end{align}
For $t\in[t_n,t_{n+1}]$,
$
(\Ic\psi)(t)=\psi_n+(t-t_n)(\delta\psi)_n$ and
$
\psi_n^\theta=\psi_n+\theta h_n(\delta\psi)_n$, which implies
\[
\abs{(\Ic\psi)(t)-\psi_n}+\abs{(\Ic\psi)(t)-\psi_n^\theta}\le 2Kh_n.
\]
For \(t\in[t_n,t_{n+1}]\), both \((I_\Delta\psi)(t)\) and \(\psi_n^\theta\) are convex combinations of \(\psi_n\) and \(\psi_{n+1}\).  Accordingly, $|(I_\Delta\psi)(t)|\le R$ and \(|\psi_n^\theta|\le R\).
The bounds on $A,b,DA,Db$ in \eqref{eq:local-elliptic} and
\eqref{eq:local-coefficient-bounds} imply that $A$ and $b$ are
Lipschitz on $\abs{x}\le R$ and that
$$
\abs{r_n}+\bigl|(\delta\psi)_n-b((\Ic\psi)(t))\bigr|\le C_{R,K}.
$$
Using \eqref{eq:AF}, we estimate the difference of the two quadratic forms:
\begin{align*}
&\left|
F((\Ic\psi)(t),(\delta\psi)_n)-\tfrac12r_n^{\top}A(\psi_n)r_n
\right|\\
&\qquad\le C_{R,K}
\bigl(
\norm{A((\Ic\psi)(t))-A(\psi_n)}
+\abs{b((\Ic\psi)(t))-b(\psi_n^\theta)}
\bigr)\\
&\qquad\le C_{R,K}
\bigl(
\abs{(\Ic\psi)(t)-\psi_n}
+\abs{(\Ic\psi)(t)-\psi_n^\theta}
\bigr)
\le C_{R,K}h_n.
\end{align*}
Substituting this bound into \eqref{eq:Slcphi} and using $h_n\le\Delta$
yields
\begin{align*}
|\St(\Ic\psi)-\Sd(\psi)|
&\le C_{R,K}\sum_{n=0}^{N-1}h_n^2
\le C_{R,K}\Delta\sum_{n=0}^{N-1}h_n
=C_{R,K}T\Delta.\qedhere
\end{align*}
\end{proof}

To construct an admissible competitor from a continuous path $\phi$,
we use the sampling operator $\mathcal{P}_\Delta$ defined in \eqref{eq:P}.
For a path in $C^2([0,T];\R^d)$, the sampling error in the action
satisfies the same first-order bound. Next, we define the $C^2$ norm by 
$\norm{\phi}_{C^2}
:=\norm{\phi}_{L^\infty(0,T;\R^d)}
+\norm{\dot\phi}_{L^\infty(0,T;\R^d)}
+\norm{\ddot\phi}_{L^\infty(0,T;\R^d)}.$ 
\begin{lemma}
\label{lem:sampling}
For every $\phi\in C^2([0,T];\R^d)$ satisfying $\phi(0)=x_0$ and $\phi(T)=x_T$,
\begin{equation}
|\Sd(\Pc\phi)-\St(\phi)|\le C_{\norm{\phi}_{C^2}}\Delta.
\label{eq:sampling-consistency}
\end{equation}
\end{lemma}

\begin{proof}
Let $y:=\Ic\Pc\phi$. For $t\in(t_n,t_{n+1})$,
\[
\dot y(t)
=\frac1{h_n}\int_{t_n}^{t_{n+1}}\dot\phi(s)\ud s.
\]
Hence, we have $\norm{\dot y}_{L^\infty}\le\norm{\dot\phi}_{L^\infty}$, and
\begin{equation}
\norm{(y-\phi)|_{(t_n,t_{n+1})}}_{L^\infty}
+h_n\norm{(\dot y-\dot\phi)|_{(t_n,t_{n+1})}}_{L^\infty}
\le Ch_n^2\norm{\ddot\phi}_{L^\infty}.
\label{eq:interpolation}
\end{equation}
Since $y$ interpolates $\phi$, we also have
$\norm{y}_{L^\infty}\le\norm{\phi}_{L^\infty}$.
Since \(F\in C^1(\mathbb R^d\times\mathbb R^d;\mathbb R)\), it is Lipschitz on the compact convex set
\[
\left\{(x,v)\in\mathbb R^d\times\mathbb R^d:
|x|\le\|\phi\|_{L^\infty},\
|v|\le\|\dot\phi\|_{L^\infty}\right\}.
\]
Using \eqref{eq:interpolation} and $h_n\le T$, we have
\begin{align}
|\St(y)-\St(\phi)|
&\le\int_0^T
\bigl|F(y,\dot y)-F(\phi,\dot\phi)\bigr|\ud t
\le C_{\norm{\phi}_{C^2}}\int_0^T
\bigl(\abs{y-\phi}+\abs{\dot y-\dot\phi}\bigr)\ud t \notag\\
&\le C_{\norm{\phi}_{C^2}}\sum_{n=0}^{N-1}(h_n^3+h_n^2)
\le C_{\norm{\phi}_{C^2}}\Delta.
\label{eq:sampling-continuous-action}
\end{align}
Applying Lemma~\ref{lem:affine-consistency} with $\psi=\Pc\phi$
and using \eqref{eq:sampling-continuous-action}, we obtain
\begin{align*}
|\Sd(\Pc\phi)-\St(\phi)|
&\le |\Sd(\Pc\phi)-\St(y)|+|\St(y)-\St(\phi)|
\le C_{\norm{\phi}_{C^2}}\Delta.\qedhere
\end{align*}
\end{proof}

We next compare the discrete minimum $m_\Delta$ with the continuous
minimum $m$.

\begin{theorem}[First-order convergence of the minimum action]
\label{thm:main}
Under Assumption~\ref{ass:main}, there exists $C>0$ such that, for every
partition with $\Delta\le\Delta_0$ and every $\theta\in[0,1]$,
\begin{equation}
|m_\Delta-m|\le C\Delta.
\label{eq:main-rate}
\end{equation}
\end{theorem}

\begin{proof}
Let $\Delta\le\Delta_0$.
Proposition~\ref{prop:existence} provides continuous and discrete
minimizers $\phi^*$ and $\psi^\Delta$, with actions bounded by $M_0$.
Both minimizers are stationary. Applying the regularity estimates in
Propositions~\ref{prop:continuous-reg} and~\ref{prop:discrete-reg} leads to
\[
\norm{\phi^*}_{C^2}
+\norm{\Ic\psi^\Delta}_{W^{1,\infty}}
\le C.
\]
Then Lemma~\ref{lem:sampling} yields
\begin{equation}
m_\Delta
\le\Sd(\Pc\phi^*)
\le\St(\phi^*)+C\Delta
=m+C\Delta.
\label{eq:upper}
\end{equation}
Conversely, since $\Ic\psi^\Delta\in\Ac$,
Lemma~\ref{lem:affine-consistency} bounds the continuous minimum from above:
\begin{equation}
m
\le\St(\Ic\psi^\Delta)
\le\Sd(\psi^\Delta)+C\Delta
=m_\Delta+C\Delta.
\label{eq:lower}
\end{equation}
Combining \eqref{eq:upper} and \eqref{eq:lower} proves
\eqref{eq:main-rate}.
\end{proof}

\begin{remark}
\label{rem:inexact-minimization}
Under the hypotheses of Theorem~\ref{thm:main}, let
$\psi\in\Ac_\Delta$ be a computed path satisfying
\(
0\le\Sd(\psi)-m_\Delta\le\varepsilon_\Delta.
\)
Using \eqref{eq:main-rate}, we obtain
\[
|\Sd(\psi)-m|
\le
\underbrace{\Sd(\psi)-m_\Delta}_{\text{optimization error}}
+\underbrace{|m_\Delta-m|}_{\text{discretization error}}
\le\varepsilon_\Delta+C\Delta.
\]
In particular, if $\varepsilon_\Delta\le c\Delta$ with $c$ independent
of the partition, then
\[
|\Sd(\psi)-m|\le(C+c)\Delta.
\]
%A small gradient with respect to the interior nodes indicates
%approximate stationarity but does not, by itself, bound the global
%optimality gap $\Sd(\psi)-m_\Delta$.
\end{remark}
%\subsection{Strong convergence of minimizing paths}

Under global Lipschitz assumptions, weak subsequential convergence of
interpolated discrete minimizers in $H^1(0,T;\R^d)$ was established in
\cite[Theorem~3.3]{HongJinSheng2025}. Under Assumption~\ref{ass:main}, we can show that
interpolated discrete minimizers converge strongly along subsequences
in the same space.

\begin{proposition}
\label{prop:minimizer-convergence}
Under Assumption~\ref{ass:main}, let $\Delta_j\to0$.
For each $j$, let $\psi^{(j)}$ minimize $\mathcal S_{\Delta_j}$ over
$\Ac_{\Delta_j}$ with parameter $\theta_j\in[0,1]$, held fixed on
that partition.
Then the sequence of interpolants $(I_{\Delta_j}\psi^{(j)})_{j\ge1}$
has a subsequence satisfying
$$
I_{\Delta_{j_k}}\psi^{(j_k)}
\to \phi^\ast
\quad\text{strongly in }H^1(0,T;\R^d),
$$
where $\phi^\ast$ minimizes $\St$ over $\Ac$.
If $\St$ has a unique minimizer over $\Ac$, the entire sequence of
interpolants converges to $\phi^\ast$ strongly in $H^1(0,T;\R^d)$.
\end{proposition}

\begin{proof}
Set $y_j:=I_{\Delta_j}\psi^{(j)}$.
After discarding finitely many terms, we have $\Delta_j\le\Delta_0$. By Proposition~\ref{prop:existence},
$\mathcal S_{\Delta_j}(\psi^{(j)})\le M_0$.
Lemma~\ref{lem:discrete-sublevel} bounds $y_j$ in
$H^1(0,T;\R^d)$. Thus for some subsequence, not relabelled, we obtain
$\phi^*\in\Ac$ such that
\[
y_j\rightharpoonup\phi^*\quad\text{in }H^1(0,T;\R^d),\qquad
y_j\to\phi^*\quad\text{in }C([0,T];\R^d).
\]
For stationary discrete minimizers, Proposition~\ref{prop:discrete-reg}
bounds their difference quotients uniformly, so Lemma~\ref{lem:affine-consistency} implies
\[
|\St(y_j)-m_{\Delta_j}|\le C\Delta_j.
\]
Together with Theorem~\ref{thm:main}, this implies $\St(y_j)\to m$.
The weak lower semicontinuity argument in the proof of
Proposition~\ref{prop:existence} yields
$\St(\phi^*)\le\liminf_j\St(y_j)=m$.
As $\phi^*\in\Ac$, it is a minimizer of $\St$ over $\Ac$. Thus,
\begin{equation}
\St(y_j)\to\St(\phi^*)=m.
\label{eq:minimizer-action-convergence}
\end{equation}

To deduce strong convergence of the first derivatives from
\eqref{eq:minimizer-action-convergence}, define the associated controls by
\begin{equation}
u_j:=\sigma(y_j)^{-1}\bigl(\dot y_j-b(y_j)\bigr),
\qquad
u:=\sigma(\phi^*)^{-1}\bigl(\dot\phi^*-b(\phi^*)\bigr).
\label{eq:minimizer-controls}
\end{equation}
The ranges of $y_j$ and $\phi^*$ are contained in a common
closed ball $\mathcal K:=\overline{B_R(0)}\subset\R^d$, with $R$
independent of $j$. On this compact set, $b,\sigma,\sigma^{-1}$ are
bounded and uniformly continuous.
Since $\sup_j\norm{\dot y_j}_{L^2}<\infty$ and $y_j$ converges uniformly,
\begin{align*}
\norm{
u_j-\sigma(\phi^*)^{-1}\bigl(\dot y_j-b(\phi^*)\bigr)
}_{L^2}\le C\left(
\norm{\sigma(y_j)^{-1}-\sigma(\phi^*)^{-1}}_{L^\infty}
+\norm{b(y_j)-b(\phi^*)}_{L^\infty}
\right)\to0.
\end{align*}
Together with the weak convergence of $\dot y_j$, this implies
$u_j\rightharpoonup u$ in $L^2(0,T;\R^d)$. On the other hand, \eqref{eq:minimizer-action-convergence} gives
$$
\norm{u_j}_{L^2}^2
=2\St(y_j)
\to 2\St(\phi^*)
=\norm{u}_{L^2}^2.
$$
In the Hilbert space $L^2(0,T;\R^d)$, weak convergence together with
convergence of norms implies strong convergence. Hence
$u_j\to u$ strongly in $L^2(0,T;\R^d)$.
We now estimate the difference of the first derivatives using \eqref{eq:minimizer-controls}:
\begin{align*}
\norm{\dot y_j-\dot\phi^*}_{L^2}
&\le \sqrt T\,\norm{b(y_j)-b(\phi^*)}_{L^\infty}\\
&\quad
+\norm{\sigma(y_j)-\sigma(\phi^*)}_{L^\infty}
 \norm{u_j}_{L^2}\\
&\quad
+\norm{\sigma(\phi^*)}_{L^\infty}\norm{u_j-u}_{L^2}
\to0.
\end{align*}
Together with uniform convergence of $y_j$ to $\phi^*$, this proves
$y_j\to\phi^*$ strongly in $H^1(0,T;\R^d)$.
If $\St$ has a unique minimizer over $\Ac$, applying the same argument
to every subsequence of $(y_j)_{j\ge1}$ yields strong convergence of the
entire interpolant sequence to $\phi^*$ in $H^1(0,T;\R^d)$.
\end{proof}
\subsection{First-order bias and sharpness}
\label{sec:bias}

\leavevmode In this part, on uniform grids, we identify the leading bias in the minimum
action, derive first-order error bounds for minimizing paths, and
establish the sharpness of the minimum-action estimate.
Define
\begin{align}
G_\theta(x,v)
&:=(\tfrac12-\theta)(v-b(x))^{\top}A(x)Db(x)v
-\tfrac14(v-b(x))^{\top}DA(x)[v](v-b(x)),
\label{eq:bias-density}\\
\mathcal B_\theta(\phi)
&:=\int_0^T G_\theta(\phi,\dot\phi)\ud t,
\qquad \phi\in W^{1,\infty}(0,T;\R^d).
\label{eq:bias-functional}
\end{align}

For piecewise affine interpolants on uniform time grids, the following
lemma shows that $\mathcal B_\theta$ is the coefficient of the first-order
discrepancy between the discrete and continuous actions.

\begin{lemma}
\label{lem:bias-cell}
Assume that $b\in C^2(\R^d;\R^d)$,
$\sigma\in C^2(\R^d;\R^{d\times d})$, and
$\det\sigma(x)\ne0$ for every $x\in\R^d$.
Let $R,K>0$ and consider a uniform grid $t_n=n\Delta$,
where $\Delta=T/N$ and $N\ge1$.
Then for every $\theta\in[0,1]$ and every
$\psi=(\psi_n)_{n=0}^N\in\Ac_\Delta$ satisfying
\begin{equation}
\max_{0\le n\le N}\abs{\psi_n}\le R,
\qquad
\max_{0\le n<N}\abs{(\delta\psi)_n}\le K,
\label{eq:bias-nodal-bounds}
\end{equation}
we have
\begin{equation}
\left|
\Sd(\psi)-\St(\Ic\psi)
-\Delta\mathcal B_\theta(\Ic\psi)
\right|
\le C_{R,K}\Delta^2,
\label{eq:cell-bias}
\end{equation}
where $C_{R,K}$ is independent of $N$, $\theta$, and $\psi$.
\end{lemma}

\begin{proof}
For $0\le n<N$, set
$
r_n:=(\delta\psi)_n-b(\psi_n^\theta).
$
By \eqref{eq:discrete-action-quadratic}, \eqref{eq:action-integral}, and the definition of $\Ic$,
\begin{align}
\Sd(\psi)-\St(\Ic\psi) =
\sum_{n=0}^{N-1}
\left[
\frac{\Delta}{2}r_n^{\top}A(\psi_n)r_n
-\int_0^\Delta
F(\psi_n+\tau(\delta\psi)_n,(\delta\psi)_n)\ud\tau
\right].
\label{eq:bias-action-difference}
\end{align}
We compare the two terms in each summand with
$\Delta F(\psi_n,(\delta\psi)_n)$.

For $0\le\tau\le\Delta$, the affine interpolant satisfies
$
\psi_n+\tau(\delta\psi)_n
=(\Ic\psi)(t_n+\tau).
$
Hence, the nodal bounds in \eqref{eq:bias-nodal-bounds} imply
\begin{equation}
\abs{\psi_n+\tau(\delta\psi)_n}\le R,
\qquad
\abs{(\delta\psi)_n}\le K.
\label{eq:bias-interpolant-bounds}
\end{equation}
The $C^2$ regularity of $b,\sigma$ and the invertibility of $\sigma$
imply that $A,b$ and their derivatives through order two are bounded
on $\{x\in\R^d:\abs{x}\le R\}$, which yields
\[
\abs{\nabla_xF(x,v)}+\norm{D_x^2F(x,v)}
\le C_{R,K},
\qquad \abs{x}\le R,\ \abs{v}\le K.
\]

For the continuous action, \eqref{eq:bias-interpolant-bounds} and
the bound on $D_x^2F$ imply, for $0\le\tau\le\Delta$,
\begin{align*}
&\left|
F(\psi_n+\tau(\delta\psi)_n,(\delta\psi)_n)
-F(\psi_n,(\delta\psi)_n)
-\tau\nabla_xF(\psi_n,(\delta\psi)_n)\cdot(\delta\psi)_n
\right|\le C_{R,K}\tau^2.
\end{align*}
Integrating over $[0,\Delta]$ yields the continuous cell expansion:
\begin{align*}
\bigg|
\int_0^\Delta
F(\psi_n+\tau(\delta\psi)_n,(\delta\psi)_n)\ud\tau
-\Delta F(\psi_n,(\delta\psi)_n)
-\frac{\Delta^2}{2}
\nabla_xF(\psi_n,(\delta\psi)_n)\cdot(\delta\psi)_n
\bigg|
\le C_{R,K}\Delta^3.
\end{align*}

For the discrete action, since
$\psi_n^\theta=\psi_n+\theta\Delta(\delta\psi)_n$, it follows from
\eqref{eq:bias-interpolant-bounds} and the boundedness of $Db,D^2b$ that
\begin{gather*}
\abs{r_n-\widehat r_n}
=\abs{b(\psi_n^\theta)-b(\psi_n)}
\le C_{R,K}\Delta,\\
\abs{r_n-\widehat r_n
+\theta\Delta Db(\psi_n)(\delta\psi)_n}
\le C_{R,K}\Delta^2,
\end{gather*}
where $\widehat r_n:=(\delta\psi)_n-b(\psi_n)$. Using the identity
\begin{align*}
\frac{1}{2}r_n^{\top}A(\psi_n)r_n
= F(\psi_n,(\delta\psi)_n)
+\widehat r_n^{\top}A(\psi_n)(r_n-\widehat r_n)+\frac{1}{2}
(r_n-\widehat r_n)^{\top}A(\psi_n)(r_n-\widehat r_n),
\end{align*}
we obtain the discrete cell expansion:
\begin{align*}
&\left|
\frac{\Delta}{2}r_n^{\top}A(\psi_n)r_n
-\Delta F(\psi_n,(\delta\psi)_n)
+\theta\Delta^2\widehat r_n^{\top}A(\psi_n)
Db(\psi_n)(\delta\psi)_n
\right|\\
&\le
\Delta\abs{\widehat r_n}\norm{A(\psi_n)}
\abs{r_n-\widehat r_n
+\theta\Delta Db(\psi_n)(\delta\psi)_n}
+\frac{\Delta}{2}\norm{A(\psi_n)}
\abs{r_n-\widehat r_n}^2\\
&\le C_{R,K}\Delta^3.
\end{align*}

Substituting the preceding two estimates into
\eqref{eq:bias-action-difference}, the terms
$\Delta F(\psi_n,(\delta\psi)_n)$ cancel.
By \eqref{eq:metric-derivative} with $\xi=(\delta\psi)_n$,
the coefficient of $\Delta^2$ is
\begin{align*}
&-\theta\widehat r_n^{\top}A(\psi_n)Db(\psi_n)(\delta\psi)_n
-\tfrac12\nabla_xF(\psi_n,(\delta\psi)_n)\cdot(\delta\psi)_n\\
&\quad=(\tfrac12-\theta)\widehat r_n^{\top}A(\psi_n)
Db(\psi_n)(\delta\psi)_n
-\tfrac14\widehat r_n^{\top}DA(\psi_n)[(\delta\psi)_n]\widehat r_n\\
&\quad=G_\theta(\psi_n,(\delta\psi)_n).
\end{align*}
Accordingly, it follows that
\begin{equation}
\Big|
\Sd(\psi)-\St(\Ic\psi)
-\Delta^2\sum_{n=0}^{N-1}G_\theta(\psi_n,(\delta\psi)_n)
\Big|\le C_{R,K}N\Delta^3
=C_{R,K}T\Delta^2.
\label{eq:bias-summed-error}
\end{equation}

It remains to compare the nodal sum with the integral defining
$\mathcal B_\theta(\Ic\psi)$.
By \eqref{eq:bias-density}, the derivatives of $G_\theta$ involve
only $A,b$ and their derivatives through order two.
Hence $G_\theta$ is Lipschitz for $\abs{x}\le R$ and $\abs{v}\le K$,
uniformly in $\theta\in[0,1]$.
Using \eqref{eq:bias-interpolant-bounds}, we obtain
\begin{align*}
\bigg|
\mathcal B_\theta(\Ic\psi)
-\Delta\sum_{n=0}^{N-1}G_\theta(\psi_n,(\delta\psi)_n)
\bigg|&\le
\sum_{n=0}^{N-1}\int_0^\Delta
\left|
G_\theta(\psi_n+\tau(\delta\psi)_n,(\delta\psi)_n)
-G_\theta(\psi_n,(\delta\psi)_n)
\right|\ud\tau\\
&\le C_{R,K}\sum_{n=0}^{N-1}\int_0^\Delta\tau\ud\tau
\le C_{R,K}T\Delta.
\end{align*}
Multiplying the preceding bound by $\Delta$ and combining it with
\eqref{eq:bias-summed-error} completes the proof.
\end{proof}

In the following theorem, we identify the leading
coefficient of the minimum-action error on uniform grids.

\begin{theorem}[First-order asymptotic coefficient]
\label{thm:bias}
Under Assumption~\ref{ass:main}, suppose that
$b\in C^2(\R^d;\R^d)$,
$\sigma\in C^2(\R^d;\R^{d\times d})$, and
$\St$ has a unique minimizer $\phi^*$ over $\Ac$.
Fix \(\theta\in[0,1]\). On uniform grids \(t_n=nT/N\), with \(\Delta=T/N\), as \(N\to\infty\), we have 
\begin{equation}
m_\Delta
=m+\Delta\mathcal B_\theta(\phi^*)+o(\Delta).
\label{eq:bias-theorem}
\end{equation}
\end{theorem}

\begin{proof}
On a uniform grid with $\Delta\le\Delta_0$, let $\psi^\Delta$
minimize $\Sd$ over $\Ac_\Delta$, and set
$y_\Delta:=\Ic\psi^\Delta$.
We first derive a lower bound for $(m_\Delta-m)/\Delta$.
By \eqref{eq:comparison-level} and Proposition~\ref{prop:discrete-reg},
$\norm{y_\Delta}_{W^{1,\infty}}\le C$.
Thus \eqref{eq:bias-nodal-bounds} holds for $\psi^\Delta$ with
constants independent of $\Delta$. Applying \eqref{eq:cell-bias}
and using $\Sd(\psi^\Delta)=m_\Delta$, we obtain
\begin{equation}
\left|m_\Delta-\St(y_\Delta)-\Delta\mathcal B_\theta(y_\Delta)\right|
\le C\Delta^2.
\label{eq:bias-minimizer-expansion}
\end{equation}
Since $y_\Delta\in\Ac$, we have $\St(y_\Delta)\ge m$.
Subtracting $m$ and dividing by $\Delta$, we obtain
\begin{align}
\frac{m_\Delta-m}{\Delta}
\ge\frac{\St(y_\Delta)-m}{\Delta}
+\mathcal B_\theta(y_\Delta)-C\Delta \ge\mathcal B_\theta(y_\Delta)-C\Delta.
\label{eq:bias-lower}
\end{align}

Since $\phi^*$ is the unique minimizer of $\St$ over $\Ac$,
Proposition~\ref{prop:minimizer-convergence} implies that
$y_\Delta\to\phi^*$ in $H^1(0,T;\R^d)$.
To pass to the limit in \eqref{eq:bias-lower}, we compare
$\mathcal B_\theta(y_\Delta)$ with $\mathcal B_\theta(\phi^*)$.
By the common action bound \eqref{eq:comparison-level} and
Proposition~\ref{prop:continuous-reg}, we also have
$\norm{\phi^*}_{W^{1,\infty}}\le C$.
By \eqref{eq:bias-density} and the $C^2$ regularity of $A,b$,
$\nabla_xG_\theta$ and $\nabla_vG_\theta$ are bounded on every
product of closed balls in $\R^d\times\R^d$.
Fix $R,V>0$ and let
$\phi_1,\phi_2\in W^{1,\infty}(0,T;\R^d)$ satisfy
$\norm{\phi_i}_{L^\infty}\le R$ and
$\norm{\dot\phi_i}_{L^\infty}\le V$, $i=1,2$.
The mean value estimate and the Cauchy--Schwarz inequality give
\begin{align}
|\mathcal B_\theta(\phi_1)-\mathcal B_\theta(\phi_2)|
&=\left|\int_0^T
\bigl[G_\theta(\phi_1,\dot\phi_1)-G_\theta(\phi_2,\dot\phi_2)\bigr]
\ud t\right| \notag\\
&\le C_{R,V}\int_0^T
\bigl(\abs{\phi_1-\phi_2}+\abs{\dot\phi_1-\dot\phi_2}\bigr)\ud t \notag\\
&\le C_{R,V}\norm{\phi_1-\phi_2}_{H^1}.
\label{eq:bias-functional-comparison}
\end{align}
In particular,
\(
|\mathcal B_\theta(y_\Delta)-\mathcal B_\theta(\phi^*)|
\le C\norm{y_\Delta-\phi^*}_{H^1}\to0,
\) as $\Delta\to0$.
Thus \eqref{eq:bias-lower} yields
\begin{equation}
\liminf_{\Delta\to0}\frac{m_\Delta-m}{\Delta}
\ge\mathcal B_\theta(\phi^*).
\label{eq:bias-liminf}
\end{equation}

For an upper bound, we choose the sampled path $\Pc\phi^*\in\Ac_\Delta$
and set $z_\Delta:=\Ic\Pc\phi^*$.
Using the second derivative bound from
Proposition~\ref{prop:continuous-reg} in \eqref{eq:interpolation}
on each interval $[t_n,t_{n+1}]$, we obtain
\begin{equation}
\|z_\Delta-\phi^*\|_{H^1}\le C\Delta,\qquad
\|z_\Delta\|_{W^{1,\infty}}\le C.
\label{eq:comparison-interpolation-bounds}
\end{equation}
In particular, $\Pc\phi^*$ satisfies \eqref{eq:bias-nodal-bounds}
with constants independent of $\Delta$.
Discrete minimality and \eqref{eq:cell-bias} imply
\begin{equation}
m_\Delta\le\Sd(\Pc\phi^*)
\le\St(z_\Delta)+\Delta\mathcal B_\theta(z_\Delta)+C\Delta^2.
\label{eq:bias-competitor-comparison}
\end{equation}
To use this comparison, we estimate $\St(z_\Delta)-m$ and
$|\mathcal B_\theta(z_\Delta)-\mathcal B_\theta(\phi^*)|$.
Noting that interpolation preserves the endpoints, we have
$z_\Delta-\phi^*\in H_0^1(0,T;\R^d)$ and
$\phi^*+s(z_\Delta-\phi^*)\in\Ac$ for every $s\in\R$.
Since $\phi^*$ minimizes $\St$, differentiation at $s=0$ yields
\[
D\St(\phi^*)[z_\Delta-\phi^*]=0.
\]
The paths $z_\Delta,\phi^*$ and their first derivatives are
uniformly bounded. Hence, for almost every $t\in(0,T)$, the line segment
in $\R^d\times\R^d$ joining $(\phi^*(t),\dot\phi^*(t))$ to
$(z_\Delta(t),\dot z_\Delta(t))$ lies in the compact set
$\mathcal K:=\overline{B_R(0)}\times\overline{B_V(0)}$, with $R,V$
independent of $t$ and $\Delta$. By \eqref{eq:AF} and the $C^2$
regularity of $b,\sigma$, we have
$\sup_{(x,v)\in\mathcal K}\norm{D^2F(x,v)}<\infty$, so the Hessian
bound is uniform in $t$ and $\Delta$.
Applying Taylor's formula, one has 
\begin{align*}
&\left|\St(z_\Delta)-\St(\phi^*)
-D\St(\phi^*)[z_\Delta-\phi^*]\right|\\
&=
\bigg|\int_0^T\Bigl[
F(z_\Delta,\dot z_\Delta)-F(\phi^*,\dot\phi^*)
-\nabla_xF(\phi^*,\dot\phi^*)\cdot(z_\Delta-\phi^*)
-\nabla_vF(\phi^*,\dot\phi^*)\cdot
(\dot z_\Delta-\dot\phi^*)\Bigr]\ud t\bigg|\\
&\le C\int_0^T
\bigl(\abs{z_\Delta-\phi^*}^2+\abs{\dot z_\Delta-\dot\phi^*}^2\bigr)\ud t.
\end{align*}
Since $\St(\phi^*)=m$ and $z_\Delta\in\Ac$,
\eqref{eq:comparison-interpolation-bounds} yields
\begin{equation}
0\le\St(z_\Delta)-m
\le C\|z_\Delta-\phi^*\|_{H^1}^2
\le C\Delta^2.
\label{eq:stationary-interpolation}
\end{equation}
By \eqref{eq:bias-functional-comparison} and
\eqref{eq:comparison-interpolation-bounds},
$|\mathcal B_\theta(z_\Delta)-\mathcal B_\theta(\phi^*)|\le C\Delta$.
Combining this and \eqref{eq:stationary-interpolation}
with \eqref{eq:bias-competitor-comparison} yields
\begin{equation}
\frac{m_\Delta-m}{\Delta}
\le\frac{\St(z_\Delta)-m}{\Delta}
+\mathcal B_\theta(z_\Delta)+C\Delta
\le\mathcal B_\theta(\phi^*)+C\Delta.
\label{eq:bias-upper}
\end{equation}
Taking the upper limit in \eqref{eq:bias-upper} and using
\eqref{eq:bias-liminf} proves \eqref{eq:bias-theorem}.
\end{proof}

We now estimate the error in the minimizing paths, assuming local
quadratic growth of $\St$ at its unique minimizer.

\begin{theorem}[First-order convergence of minimizing paths]
\label{thm:path-rate}
Under Assumption~\ref{ass:main}, suppose that
$b\in C^2(\R^d;\R^d)$,
$\sigma\in C^2(\R^d;\R^{d\times d})$, and
$\St$ has a unique minimizer $\phi^*$ over $\Ac$.
Suppose further that, for some $c_*,\varepsilon_*>0$,
\begin{equation}
\St(\phi)-m\ge c_*\norm{\phi-\phi^*}_{H^1}^2
\quad\text{for all }\phi\in\Ac
\text{ with }\norm{\phi-\phi^*}_{H^1}\le\varepsilon_*.
\label{eq:path-quadratic-growth}
\end{equation}
Then there exist $C>0$ and $\Delta_1>0$ such that every global
minimizer $\psi^\Delta$ of $\Sd$ over $\Ac_\Delta$
on a uniform grid with
$\Delta\le\Delta_1$ satisfies
\begin{equation}
\norm{\Ic\psi^\Delta-\phi^*}_{H^1}
+\norm{\Ic\psi^\Delta-\phi^*}_{L^\infty}
\le C\Delta.
\label{eq:path-first-order}
\end{equation}
\end{theorem}

\begin{proof}
Let $\psi^\Delta$ minimize $\Sd$ over $\Ac_\Delta$ on a uniform grid
with $\Delta\le\Delta_0$, and set $y_\Delta:=\Ic\psi^\Delta\in\Ac$.
Since $\phi^*$ is the unique minimizer of $\St$ over $\Ac$,
Proposition~\ref{prop:minimizer-convergence} implies that
$\norm{y_\Delta-\phi^*}_{H^1}\to0$ as $\Delta\to 0$.
Consequently, we can choose $0<\Delta_1\le\Delta_0$ so that
$\norm{y_\Delta-\phi^*}_{H^1}\le\varepsilon_*$ for all $\Delta\le\Delta_1$ and all
choices of $\psi^\Delta$. Then \eqref{eq:path-quadratic-growth} gives
\[
c_*\norm{y_\Delta-\phi^*}_{H^1}^2\le\St(y_\Delta)-m.
\]

To estimate $\St(y_\Delta)-m$ from above, we use the admissible competitor
$\Pc\phi^*\in\Ac_\Delta$ and set $z_\Delta:=\Ic\Pc\phi^*$.
It follows from comparisons \eqref{eq:bias-minimizer-expansion} and
\eqref{eq:bias-competitor-comparison} that
\begin{align}
\St(y_\Delta)+\Delta\mathcal B_\theta(y_\Delta)
\le m_\Delta+C\Delta^2 \le\St(z_\Delta)+\Delta\mathcal B_\theta(z_\Delta)+C\Delta^2.
\label{eq:path-discrete-comparison}
\end{align}
To control the bias difference in \eqref{eq:path-discrete-comparison},
note that \eqref{eq:comparison-level},
Proposition~\ref{prop:discrete-reg}, and
\eqref{eq:comparison-interpolation-bounds} ensure
$\norm{y_\Delta}_{W^{1,\infty}}+\norm{z_\Delta}_{W^{1,\infty}}\le C$.
Applying \eqref{eq:bias-functional-comparison} and then
\eqref{eq:comparison-interpolation-bounds}, we obtain
\begin{align*}
|\mathcal B_\theta(z_\Delta)-\mathcal B_\theta(y_\Delta)|
\le C\norm{z_\Delta-y_\Delta}_{H^1}
%&\le C\bigl(\norm{z_\Delta-\phi^*}_{H^1}
%+\norm{y_\Delta-\phi^*}_{H^1}\bigr)\\
\le C\bigl(\Delta+\norm{y_\Delta-\phi^*}_{H^1}\bigr).
\end{align*}
Subtracting $m$ from \eqref{eq:path-discrete-comparison} and using
\eqref{eq:stationary-interpolation} with this bound on the bias
difference, we obtain
\begin{align}
0\le\St(y_\Delta)-m
&\le\St(z_\Delta)-m
+\Delta\bigl(\mathcal B_\theta(z_\Delta)
-\mathcal B_\theta(y_\Delta)\bigr)+C\Delta^2 \notag\\
&\le C\Delta^2+C\Delta\norm{y_\Delta-\phi^*}_{H^1}.
\label{eq:path-bias-comparison}
\end{align}

Combining \eqref{eq:path-quadratic-growth} with
\eqref{eq:path-bias-comparison} and applying Young's inequality yields
\begin{align*}
c_*\norm{y_\Delta-\phi^*}_{H^1}^2\le \St(y_\Delta)-m
\le C\Delta^2+C\Delta\norm{y_\Delta-\phi^*}_{H^1}\le \frac{c_*}{2}\norm{y_\Delta-\phi^*}_{H^1}^2+C\Delta^2.
\end{align*}
Rearranging this inequality, we have $\norm{y_\Delta-\phi^*}_{H^1}\le C\Delta$.
Since $y_\Delta(0)=\phi^*(0)$,
\[
\norm{y_\Delta-\phi^*}_{L^\infty}
\le\sqrt T\,\norm{\dot y_\Delta-\dot\phi^*}_{L^2}
\le C\Delta.\qedhere
\]
\end{proof}

The following two cases ensure uniqueness of the minimizer and
the local quadratic growth condition \eqref{eq:path-quadratic-growth}
in Theorem~\ref{thm:path-rate}.

\begin{remark}
\label{rem:path-stability}
Under Assumption~\ref{ass:main} with $b$ and $\sigma$ of class $C^2$, let
$\phi^*$ be a minimizer of $\St$ over $\Ac$.
The embedding $H^1(0,T;\R^d)\hookrightarrow C([0,T];\R^d)$ and
the quadratic dependence of $F$ on its second argument imply
$\St\in C^2(H^1(0,T;\R^d);\R)$.
By Taylor's formula and the stationarity of $\phi^*$, the bound
\begin{equation}
D^2\St(\phi^*)[\eta,\eta]
\ge\mu\norm{\eta}_{H^1}^2
\quad\text{for all }\eta\in H_0^1(0,T;\R^d),
\qquad \mu>0
\label{eq:path-second-variation}
\end{equation}
implies \eqref{eq:path-quadratic-growth}.

(1) Suppose that $d=1$ and $b=0$. Define
\(
H(x):=\int_{x_0}^x\frac{1}{|\sigma(r)|}\ud r.
\)
Since $\sigma$ is continuous and nowhere zero, $H$ is strictly
increasing, and
\begin{equation}
\St(\phi)=\frac12\int_0^T|(H\circ\phi)'|^2\ud t.
\label{eq:scalar-transformed-action}
\end{equation}
The Cauchy--Schwarz inequality identifies the unique minimizer given by $
\phi^*(t)=H^{-1}\left(\frac{t}{T}H(x_T)\right)$, with $
(H\circ\phi^*)'=\frac{H(x_T)}{T}.
$
For $\eta\in H_0^1(0,T;\R)$, the first two derivatives of
$H(\phi^*+\varepsilon\eta)$ with respect to $\varepsilon$ at zero
are $H'(\phi^*)\eta$ and
$H''(\phi^*)\eta^2$. Differentiating the action in
\eqref{eq:scalar-transformed-action} twice with respect to
$\varepsilon$, we obtain
\begin{align*}
D^2\St(\phi^*)[\eta,\eta]
&=\Big.\frac{\mathrm d^2}{\mathrm d\varepsilon^2}
\St(\phi^*+\varepsilon\eta)\Big|_{\varepsilon=0}\\
&=\norm{(H'(\phi^*)\eta)'}_{L^2}^2
+\int_0^T(H\circ\phi^*)'
\bigl(H''(\phi^*)\eta^2\bigr)'\ud t\\
&=\norm{(H'(\phi^*)\eta)'}_{L^2}^2
+\frac{H(x_T)}{T}\int_0^T
\bigl(H''(\phi^*)\eta^2\bigr)'\ud t\\
%&=\norm{(H'(\phi^*)\eta)'}_{L^2}^2
%+\frac{H(x_T)}{T}\bigl[H''(\phi^*)\eta^2\bigr]_0^T\\
&=\norm{(H'(\phi^*)\eta)'}_{L^2}^2.
\end{align*}
The last equality uses the fundamental theorem of calculus and
$\eta(0)=\eta(T)=0$.
Since $\phi^*\in C^1([0,T];\R)$ and $H'>0$, the functions
$H'(\phi^*)$ and $1/H'(\phi^*)$ belong to $W^{1,\infty}(0,T;\R)$.
Multiplication by $1/H'(\phi^*)$ and Poincar\'e's inequality yield \eqref{eq:path-second-variation}, namely, 
\[
\norm{\eta}_{H^1}
\le C\norm{H'(\phi^*)\eta}_{H^1}
\le C\norm{(H'(\phi^*)\eta)'}_{L^2}
=C\bigl(D^2\St(\phi^*)[\eta,\eta]\bigr)^{1/2}.
\]
%This proves \eqref{eq:path-second-variation}.

(2) Suppose that $\sigma(x)\equiv\sigma_0$, where
$\sigma_0\in\R^{d\times d}$ is invertible, and
$b(x)=Bx+c$ for a constant matrix $B\in\R^{d\times d}$ and
a vector $c\in\R^d$. In this case, the action is quadratic with
\[
D^2\St(\phi)[\eta,\eta]
=\norm{\sigma_0^{-1}(\dot\eta-B\eta)}_{L^2}^2,
\qquad
\phi\in\Ac,\quad \eta\in H_0^1(0,T;\R^d).
\]
Since $\eta(0)=0$, variation of constants yields
\[
\eta(t)=\int_0^t e^{B(t-s)}(\dot\eta(s)-B\eta(s))\ud s.
\]
The bound $\norm{e^{B(t-s)}}\le e^{T\norm{B}}$ and the
Cauchy--Schwarz inequality imply, for $0\le t\le T$,
\begin{align*}
|\eta(t)|\le e^{T\norm{B}}\int_0^t|\dot\eta(s)-B\eta(s)|\ud s\le\sqrt T\,e^{T\norm{B}}\norm{\dot\eta-B\eta}_{L^2}.
\end{align*}
Taking the $L^2$ norm with respect to $t$, we obtain
\(
\norm{\eta}_{L^2}
\le T e^{T\norm{B}}\norm{\dot\eta-B\eta}_{L^2}.
\)
Since $\dot\eta=(\dot\eta-B\eta)+B\eta$, it follows that
\begin{align*}
\norm{\dot\eta}_{L^2}
\le\norm{\dot\eta-B\eta}_{L^2}
+\norm{B}\norm{\eta}_{L^2}\le\bigl(1+T\norm{B}e^{T\norm{B}}\bigr)
\norm{\dot\eta-B\eta}_{L^2}.
\end{align*}
Combining these estimates yields
\[
\norm{\eta}_{H^1}
\le C_{B,T}\norm{\dot\eta-B\eta}_{L^2}
\le C_{B,T}\norm{\sigma_0}
\norm{\sigma_0^{-1}(\dot\eta-B\eta)}_{L^2}.
\]
Thus \eqref{eq:path-second-variation} holds at every $\phi\in\Ac$.
Proposition~\ref{prop:existence} guarantees the existence of a
minimizer of $\St$ over $\Ac$. In addition,
the strict convexity of $\St$ on $\Ac$ makes the minimizer unique.

\end{remark}

The following scalar diffusion has a nonzero first-order bias and shows
that the first-order estimate for the minimum action is sharp.

\begin{example}
\label{ex:sharp}
Let $d=1$, $b=0$, $\sigma(x)=(1+x^2)^{-1/2}$,
$x_0=0$, and $x_T=1$. Define
\[
H(x):=\int_0^x\sqrt{1+s^2}\ud s
=\tfrac12\bigl(x\sqrt{1+x^2}+\operatorname{arsinh}x\bigr).
\]
For $\phi\in\Ac$, set $Y:=H\circ\phi$.
Since $Y(0)=0$ and $Y(T)=H(1)$, the Cauchy--Schwarz inequality gives
\[
\St(\phi)
=\tfrac12\int_0^T|\dot Y|^2\ud t
\ge\frac{|Y(T)-Y(0)|^2}{2T}
=\frac{H(1)^2}{2T}.
\]
Equality holds precisely when $Y(t)=H(1)t/T$.
Thus the minimizer is $\phi^*(t)=H^{-1}\bigl(H(1)t/T\bigr)$,
and the corresponding minimum action is $m=\frac{H(1)^2}{2T}$.
Then it holds that
\begin{align*}
\dot\phi^*(t)
=\frac{H(1)}{T\sqrt{1+\phi^*(t)^2}}.
\end{align*}
Since $A(x)=1+x^2$ and $b=0$, the bias $\mathcal B_\theta(\phi^*)$
can be evaluated explicitly using the change of variables
$x=\phi^*(t)$:
\begin{align*}
\mathcal B_\theta(\phi^*)
=-\tfrac14\int_0^T
2\phi^*(t)\dot\phi^*(t)^3\ud t=-\frac{H(1)^2}{4T^2}
\int_0^1\frac{2x}{1+x^2}\ud x=-\frac{H(1)^2\log2}{4T^2}<0.
\end{align*}
Therefore, Theorem~\ref{thm:bias}  shows that the first-order rate
is sharp for every $\theta\in[0,1]$. Besides,
$b$ and $\sigma$ satisfy the conditions in
Remark~\ref{rem:path-stability}(1). Hence
\eqref{eq:path-quadratic-growth} holds, and
Theorem~\ref{thm:path-rate} applies to this example.
\end{example}
\begin{corollary}[Second-order midpoint approximation of the minimum action]
\label{cor:midpoint-minimum-action}
Under Assumption~\ref{ass:main}, suppose that $b$ and $\sigma$ are of
class $C^2$ and that $A$ is constant. For $\theta=1/2$, there exists a
constant $C>0$ such that on every uniform grid with $\Delta\le\Delta_0$,
\[
|m_\Delta-m|\le C\Delta^2.
\]
In particular, this estimate applies to additive noise.
\end{corollary}

\begin{proof}
Since $DA=0$ and $\theta=1/2$, \eqref{eq:bias-density} and
\eqref{eq:bias-functional} give $\mathcal B_{1/2}\equiv0$.
Let $\psi^\Delta$ be a discrete minimizer and $\phi^*$ a continuous
minimizer. By Propositions~\ref{prop:existence}, \ref{prop:continuous-reg},
and~\ref{prop:discrete-reg}, both $\psi^\Delta$ and $\Pc\phi^*$ satisfy
\eqref{eq:bias-nodal-bounds} uniformly in $\Delta$.
Applying Lemma~\ref{lem:bias-cell} to $\psi^\Delta$ gives
\[
m_\Delta\ge\St(\Ic\psi^\Delta)-C\Delta^2\ge m-C\Delta^2.
\]
Applying the same lemma to $\Pc\phi^*$ and using the interpolation
estimate \eqref{eq:stationary-interpolation}, we obtain
\[
m_\Delta\le\Sd(\Pc\phi^*)
\le\St(\Ic\Pc\phi^*)+C\Delta^2\le m+C\Delta^2. \qedhere
\]
\end{proof}

% ===== Section 5: Numerical experiment =====
\section{Numerical experiment}
\label{sec:numerics}

In this section, we consider the two-protein genetic switch studied in
 \cite{RomaEtAl2005}. We adopt the reaction network and parameter values in
 \cite[Sec.~IV.B]{GrafkeVandenEijnden2019} and construct its chemical
 Langevin approximation.
 To obtain a $C^2$-regularized diffusion of the form \eqref{eq:sde}
 with $d=2$ satisfying our assumptions on $\R^2$, we introduce the
 following extension and fix $\delta=10^{-3}$ throughout the experiment:
\begin{equation}
\chi_\delta(s):=
\begin{cases}
s,&s\ge\delta,\\[2mm]
\displaystyle
\delta+\frac{s-\delta}
{\sqrt{1+4(s-\delta)^2/\delta^2}},&s<\delta.
\end{cases}
\label{eq:toggle-extension-map}
\end{equation}
The production rates are
\begin{equation}
f_1(x)=\frac{156}{1+\chi_\delta(x_2)^3},
\qquad
f_2(x)=\frac{30}{1+\chi_\delta(x_1)},
\qquad x\in\R^2,
\label{eq:toggle-extended-production}
\end{equation}
and the degradation rate of protein $i$ is $\chi_\delta(x_i)$.
Set
\begin{gather}
b_i(x)=f_i(x)-\chi_\delta(x_i),
\qquad
a_i(x)=f_i(x)+\chi_\delta(x_i),
\qquad i=1,2,
\label{eq:toggle-extended-coefficients}\\
\sigma(x)=\operatorname{diag}
\bigl(\sqrt{a_1(x)},\sqrt{a_2(x)}\bigr).
\label{eq:toggle-rates}
\end{gather}

The variables $x_1,x_2$ are dimensionless protein concentrations.
For $x_1,x_2\ge\delta$, the production rates in
\eqref{eq:toggle-extended-production} reduce to
$156/(1+x_2^3)$ and $30/(1+x_1)$, respectively, and the degradation
rates are $x_1,x_2$, so the regularization preserves the original
 reaction rates in this region.
Each protein suppresses the production of the other, giving two
distinct stable expression states.
The two branches of \eqref{eq:toggle-extension-map} agree through
their second derivatives at $\delta$, so $\chi_\delta\in C^2(\R)$.
Moreover, $\chi_\delta\ge\delta/2$ implies $a_i\ge\delta/2$. It can be verified that
the drift $b$ and diffusion matrix $\sigma$ are of class $C^2$
and have at most linear growth, and $\sigma(x)$ is invertible
for every $x\in\R^2$. Hence the regularized coefficients $b$ and $\sigma$ in \eqref{eq:toggle-extended-coefficients} and \eqref{eq:toggle-rates} satisfy
Assumption~\ref{ass:main}.

In this experiment, we take $T=8$ and consider transitions between the attracting
equilibria
$$
x^A\simeq(154.896,0.19244),
\qquad
x^B\simeq(0.00588,29.825).
$$
Protein 1 dominates at $x^A$, and protein 2 dominates at $x^B$.
The saddle between them is $S\simeq(12.171,2.278)$.
We use both endpoint pairs
$(x_0,x_T)=(x^A,x^B)$ and $(x_0,x_T)=(x^B,x^A)$
to examine the action and path errors as the time step decreases.

For each endpoint pair, we minimize the continuous action
\eqref{eq:continuous-action} over $\Ac$.
We use the discrete action \eqref{eq:intro-discrete}
with $\theta=0$ on uniform grids $t_n=n\Delta$, where
$\Delta=T/N$ and $N=40,80,160,320,640$.
For each $N$, we approximately minimize $\Sd$ over $\Ac_\Delta$
by a damped Newton method
with backtracking, using a Gauss--Newton direction when the exact
Hessian is not positive definite
\cite{NocedalWright2006}.
The initial guesses are the straight segment and two piecewise
linear paths through $S$ at times $0.35T$ and $0.65T$.
These paths serve only as initial guesses.
On each refined grid, we also start from the interpolated
coarser-grid path and retain the result with the smallest computed
action. For every retained path, the largest absolute component
of the interior-node gradient is below $2\times10^{-9}$.

Note that every minimizer $x$ of $\St$ over $\Ac$ satisfies the
Euler--Lagrange equation \eqref{eq:continuous-EL}, which, with
$p=\nabla_vF(x,\dot x)$, takes the following Hamiltonian form:
\begin{align*}
\dot x&=\nabla_pH(x,p),
&\dot p&=-\nabla_xH(x,p),
\end{align*}
where $x(0)=x_0,x(T)=x_T,$ and the Hamiltonian $H$ is given by
\[
H(x,p)=b(x)\cdot p+\frac12\sum_{i=1}^2a_i(x)p_i^2,
\qquad
p_i=\frac{\dot x_i-b_i(x)}{a_i(x)}.
\]
For each transition direction, we compute a continuous reference
path $\phi_{\rm ref}$ by adaptive collocation for this boundary-value problem.
We initialize the collocation solver from the discrete path computed
with $N=640$, retaining the values
$\psi_0,\psi_4,\ldots,\psi_{640}$.
The corresponding $161$ time points form the initial mesh, and
the retained concentrations provide the initial guess for $x$.
The initial guess for $p$ is computed from the difference quotients
of these values using the formula for $p_i$, and then interpolated
onto the initial mesh.
Let $(x_h,p_h)$ denote the collocation approximation to $(x,p)$ and set
$\phi_{\rm ref}:=x_h$. We use numerical quadrature to evaluate the reference action
\[
S_{\rm ref}:=\frac12\int_0^T\sum_{i=1}^2
a_i(x_h(t))p_{h,i}(t)^2\ud t \approx \St(x_h).
\]
The residual $\dot x_{h,i}-b_i(x_h)-a_i(x_h)p_{h,i}$ can be nonzero
between collocation points.
The reference actions are $5.74933774$ for $A\to B$ and
$10.91543440$ for $B\to A$.
Tightening the collocation tolerance from \(10^{-6}\) to \(2\times10^{-9}\) changes the reference actions by less than \(1.5\times10^{-8}\). All computed paths remain in \(x_1,x_2>\delta\).

\begin{figure}[htbp]
\centering
\includegraphics[width=.99\linewidth]{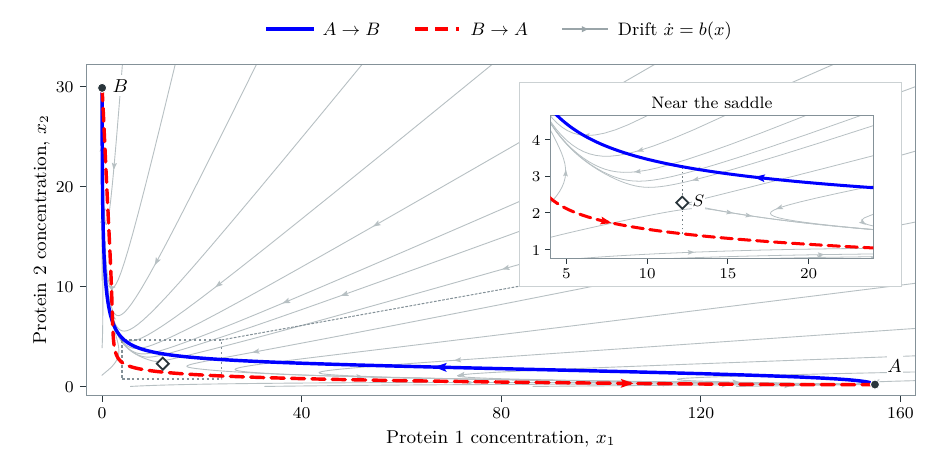}
\caption{Computed switching paths for $T=8$ on a uniform grid with
$N=640$ cells.}
\label{fig:toggle-paths}
\end{figure}

\begin{table}[H]
\centering\small
\begin{tabular}{r@{\qquad}rrrr@{\qquad}rrrr}
\toprule
 &\multicolumn{4}{c}{$A\to B$}&\multicolumn{4}{c}{$B\to A$}\\
\cmidrule(lr){2-5}\cmidrule(lr){6-9}
$N$ & $E_S$ & Rate & $E_{H^1}$ & Rate & $E_S$ & Rate & $E_{H^1}$ & Rate\\
\midrule
$40$ & $0.5756$ & -- & $21.645$ & -- & $0.0522$ & -- & $25.191$ & --\\
$80$ & $0.2963$ & $0.96$ & $11.075$ & $0.97$ & $0.0189$ & $1.47$ & $12.433$ & $1.02$\\
$160$ & $0.1504$ & $0.98$ & $5.603$ & $0.98$ & $0.0081$ & $1.23$ & $6.177$ & $1.01$\\
$320$ & $0.0758$ & $0.99$ & $2.818$ & $0.99$ & $0.0037$ & $1.12$ & $3.079$ & $1.00$\\
$640$ & $0.0380$ & $0.99$ & $1.413$ & $1.00$ & $0.0018$ & $1.06$ & $1.538$ & $1.00$\\
\bottomrule
\end{tabular}
\caption{Fixed-time action and path errors for both switching directions.
Each rate is computed as $\log_2(E_{N/2}/E_N)$ for the error column
immediately to its left.}
\label{tab:toggle-convergence}
\end{table}
For each computed discrete path $\psi$, we measure
\[
E_S=\bigl|\Sd(\psi)-S_{\rm ref}\bigr|,
\qquad
E_{H^1}=\norm{\Ic\psi-\phi_{\rm ref}}_{H^1}.
\]
Table~\ref{tab:toggle-convergence} lists the errors and the observed
rates $\log_2(E_{N/2}/E_N)$.
Both the action and path errors exhibit nearly first-order
behavior on the finest grids in both transition directions.
The comparison uses computed discrete paths and stationary
reference paths. 
%Global minimality and the uniqueness and local
%quadratic growth conditions in Theorem~\ref{thm:path-rate}
%are not verified in this experiment.
Figure~\ref{fig:toggle-paths} shows the two computed paths for
$N=640$ in the concentration plane.
The blue path goes from $A$ to $B$, and the red dashed path goes
from $B$ to $A$.
Gray arrows show the deterministic drift $\dot x=b(x)$, while the
colored arrows show the direction of travel along each path.
We use filled circles to mark the attracting equilibria $A$ and $B$, and
an open diamond to mark the saddle $S$.
The inset enlarges the boxed region on linear axes.
The paths pass on opposite sides of the saddle:
at $x_1\simeq12.17$, their second coordinates are approximately
$3.26$ and $1.43$, respectively, whereas the saddle has
$x_2\simeq2.28$.
It is observed that the two computed switching paths are not reversals of one another,
and the $B\to A$ reference action is larger than the $A\to B$
reference action.

%The companion files contain the solver, plotting script,
%optimization records, and discrete paths.

% ===== References =====

\end{document}